# MODERATE DEVIATIONS AND LAWS OF THE ITERATED LOGARITHM FOR THE LOCAL TIMES OF ADDITIVE LÉVY PROCESSES AND ADDITIVE RANDOM WALKS

By Xia Chen[1]

*University of Tennessee*

We study the upper tail behaviors of the local times of the additive Lévy processes and additive random walks. The limit forms we establish are the moderate deviations and the laws of the iterated logarithm for the $L_2$-norms of the local times and for the local times at a fixed site.

**1. Introduction.** Let $X(t)$ be a $d$-dimensional symmetric Lévy process with the characteristic exponent $\psi(\lambda)$, that is,

$$\mathbb{E}e^{i\lambda\cdot X(t)} = e^{-t\psi(\lambda)}, \qquad t \ge 0, \lambda \in \mathbb{R}^d.$$

The symmetry assumption implies that $\psi(\lambda)$ takes only real values and $\psi(\lambda) \ge 0$. Throughout we assume that there is a deterministic and positive function $a(t)$ on $\mathbb{R}^+$ such that $a(t) \to \infty$ as $t \to \infty$, and that the limit

$$\lim_{t\to\infty} t\psi\left(\frac{\lambda}{a(t)}\right) = \Psi(\lambda) \tag{1.1}$$

exists for every $\lambda \in \mathbb{R}^d$. Notice that (1.1) holds if and only if there is a symmetric $d$-dimensional stable random variable $Y$ such that

$$X(t)/a(t) \xrightarrow{d} Y. \tag{1.2}$$

In this case we have that

$$\mathbb{E}e^{i\lambda\cdot Y} = e^{-\Psi(\lambda)}, \qquad \lambda \in \mathbb{R}^d, \tag{1.3}$$

and that $\Psi(\lambda)$ is continuous, nonnegative with the properties

$$\Psi(r\lambda) = r^\alpha\Psi(\lambda) \quad \text{and} \quad \Psi(-\lambda) = \Psi(\lambda), \qquad r > 0, \lambda \in \mathbb{R}^d,$$

---

Received September 2005; revised January 2006.
[1]Supported in part by NSF Grant DMS-04-05188.

where $\alpha \in (0, 2]$ is the stable index of $Y$. It is a classic fact that $a(t)$ must be regularly varying at $\infty$ with regular exponent $1/\alpha$:

$$\lim_{t\to\infty} a(\theta t)/a(t) = \theta^{1/\alpha}, \qquad \theta > 0. \tag{1.4}$$

Throughout we assume that $Y$ is nondegenerated, which means that there is a constant $C > 0$ such that

$$C^{-1}|\lambda|^\alpha \le \Psi(\lambda) \le C|\lambda|^\alpha.$$

Let $p \ge 1$ be a fixed integer such that

$$d < \alpha p \tag{1.5}$$

and let $X_1(t), \dots, X_p(t)$ be $p$ independent copies of $X(t)$. Our first goal is to study the upper tail behaviors of the local times of the additive Lévy process

$$\overline{X}(t_1, \dots, t_p) = X_1(t_1) + \cdots + X_p(t_p), \qquad t_1, \dots, t_p \in \mathbb{R}^+.$$

The multi-parameter processes are a natural extension of existing one-parameter processes. The subject connects to other disciplines, such as functional analysis, group theory and analytic number theory; we refer to the book by Khoshnevisan [31] for these links. The multi-parameter processes also arise in applied contexts, such as mathematical statistics [40], statistical mechanics [36] and brain imaging [5].

Since they locally resemble Lévy sheets, and since they are more amenable to analysis, additive Lévy processes first arose to simplify the study of Lévy sheets (see [15, 16, 26, 27]). They also arise in the theory of intersection of Lévy trajectories. The study of additive processes connects to probabilistic potential theory. We mention [23, 25, 30, 32, 33] and refer the reader to the detailed discussion and for the further reference.

The local times of the multi-parameter process $\overline{X}(t_1, \dots, t_p)$ are defined differently in the following two different situations. The first is when

$$\int_{\mathbb{R}^d} [1 + \psi(\lambda)]^{-p}\, d\lambda < \infty. \tag{1.6}$$

According to the recent papers by Khoshnevisan, Xiao and Zhong [34, 35], the occupation measure

$$\mu_t(A) = \int_{[0,t]^p} \mathbb{1}_A(\overline{X}(s_1, \dots, s_p))\, ds_1 \cdots ds_p, \qquad A \subset \mathbb{R}^d,$$

is absolutely continuous (for all $t > 0$) with respect to the Lebesgue measure on $\mathbb{R}^d$, if and only if (1.6) holds. In this situation the local time of $\overline{X}(t_1, \dots, t_p)$, denoted by $L(t, x)$, is defined as the density function (with respect to the Lebesgue measure) of the occupation measure $\mu_t$. A formal way of writing $L(t,x)$ is

$$L(t,x)=\int_{[0,t]^p}\delta_x(X_1(s_1)+\cdots+X_p(s_p))\,ds_1\cdots ds_p,\qquad t\ge 0,\ x\in\mathbb{R}^d.$$

See [34, 35] for some discussion on the path continuity of $L(t,x)$.

In the second situation, we assume that $\psi(\lambda)$ is a periodic function with period $2\pi$:

$$\psi(\lambda+\mathbf{k}(2\pi))=\psi(\lambda),\qquad \mathbf{k}\in\mathbb{Z}^d,\ \lambda\in\mathbb{R}^d. \tag{1.7}$$

Clearly, the conditions (1.6) and (1.7) do not co-exist. Indeed, (1.7) holds if and only if the process $X(t)$ takes only $\mathbb{Z}^d$-values. Under (1.7) the local time $L(t,x)$ of $\overline{X}(t_1,\dots,t_p)$ is defined as

$$L(t,x)=\int_{[0,t]^p}\mathbb{1}_{\{X_1(s_1)+\cdots+X_p(s_p)=x\}}\,ds_1\cdots ds_p,\qquad t\ge 0,\ x\in\mathbb{Z}^d.$$

Define

$$\rho_1=\sup_{\|f\|_2=1}\int_{\mathbb{R}^d}\left[\int_{\mathbb{R}^d}\frac{f(\lambda+\gamma)f(\gamma)}{\sqrt{1+\Psi(\lambda+\gamma)}\sqrt{1+\Psi(\gamma)}}\,d\gamma\right]^p d\lambda, \tag{1.8}$$

$$\rho_2=\sup_{\|f\|_2=1}\int_{\mathbb{R}^d}\left[\int_{\mathbb{R}^d}\frac{f(\lambda+\gamma)f(\gamma)}{\sqrt{1+\Psi(\lambda+\gamma)}\sqrt{1+\Psi(\gamma)}}\,d\gamma\right]^{2p} d\lambda, \tag{1.9}$$

where

$$\|f\|_2=\left(\int_{\mathbb{R}^d}f^2(\lambda)\,d\lambda\right)^{1/2}.$$

In [8] we show that (1.5) implies that $0<\rho_1,\rho_2<\infty$.

A special case in the category defined by (1.6) is when the Lévy process is actually a stable process, in which case we always use $Y(t)$ for the stable process that generates the additive stable process $\overline{Y}$, and $L_Y(t,x)$ for the local time of $\overline{Y}$. For the process $Y(t)$, (1.2) and (1.3) automatically hold, (1.5) is equivalent to (1.6), and $Y(1)\overset{d}{=}Y$. We [9, 10] recently proved that

$$\begin{aligned}&\lim_{t\to\infty}t^{-\alpha/d}\log\mathbb{P}\{L_Y(1,0)\ge t\}\\&\quad=-\frac{d}{\alpha}(2\pi)^\alpha\left(1-\frac{d}{\alpha p}\right)^{(\alpha p-d)/d}\rho_1^{-\alpha/d}\end{aligned} \tag{1.10}$$

and that

$$\begin{aligned}&\lim_{t\to\infty}t^{-\alpha/d}\log\mathbb{P}\left\{\int_{\mathbb{R}^d}L_Y^2(1,x)\,dx\ge t\right\}\\&\quad=-\frac{d}{2\alpha}(2\pi)^\alpha\left(1-\frac{d}{2\alpha p}\right)^{(2\alpha p-d)/d}\rho_2^{-\alpha/d}.\end{aligned} \tag{1.11}$$

In the present work we consider the setting of Lévy processes, where the scaling properties given in (1.16) below are no longer reality in general. The form of large deviation we shall establish is called moderate deviation in literature, which is related to the weak law given in (1.2).

To this end, let $b_t$ be a deterministic positive function on $\mathbb{R}^+$ satisfying

$$b_t \longrightarrow \infty \quad \text{and} \quad b_t = o(t) \qquad (t \to \infty). \tag{1.12}$$

THEOREM 1.1. *Assume (1.1), (1.5) and (1.12). Under both (1.6) and (1.7),*

$$\begin{aligned}&\lim_{t\to\infty} \frac{1}{b_t} \log \mathbb{P}\left\{ L(t,0) \geq \lambda t^p a\left(\frac{t}{b_t}\right)^{-d}\right\}\\ &\quad = -(2\pi)^\alpha \frac{d}{\alpha}\left(1 - \frac{d}{\alpha p}\right)^{(\alpha p - d)/d} \left(\frac{\lambda}{\rho_1}\right)^{\alpha/d},\end{aligned} \tag{1.13}$$

$$\begin{aligned}&\lim_{t\to\infty} \frac{1}{b_t} \log \mathbb{P}\left\{ I_t \geq \lambda t^{2p} a\left(\frac{t}{b_t}\right)^{-d}\right\}\\ &\quad = -(2\pi)^\alpha \frac{d}{2\alpha}\left(1 - \frac{d}{2\alpha p}\right)^{(2\alpha p - d)/d} \left(\frac{\lambda}{\rho_2}\right)^{\alpha/d}\end{aligned} \tag{1.14}$$

*for every $\lambda > 0$, where*

$$I_t = \begin{cases} \displaystyle\int_{\mathbb{R}^d} L^2(t,x)\,dx, & \text{under (1.6)},\\ \displaystyle\sum_{x\in\mathbb{Z}^d} L^2(t,x), & \text{under (1.7)}.\end{cases} \tag{1.15}$$

Notice that in the stable case, $a(t) = t^{1/\alpha}$ and for any $t > 0$,

$$\begin{aligned} L_Y(t,0) &\stackrel{d}{=} t^{(\alpha p - d)/\alpha} L_Y(1,0) \quad \text{and}\\ \int_{\mathbb{R}^d} L_Y^2(t,x)\,dx &\stackrel{d}{=} t^{(2\alpha p - d)/\alpha} \int_{\mathbb{R}^d} L_Y^2(1,x)\,dx.\end{aligned} \tag{1.16}$$

Hence, (1.13) and (1.14) lead to (1.10) and (1.11), respectively.

Our second goal is to study the upper tail behaviors of the local times of the additive random walks. Let $S(n)$ be a symmetric $d$-dimensional random walk taking $\mathbb{Z}^d$-values with $S(0) = 0$. We assume that $S(n)$ is in the domain of attraction of a nondegenerated $d$-dimensional stable law $Y$ with the characteristic exponent $\Psi(\lambda)$. More precisely, there exist a nondecreasing deterministic positive sequence $a(n)$ with $a(n) \to \infty$ as $n \to \infty$, and a nondegenerate, symmetric $d$-dimensional stable random variable $Y$ described as

above, such that

$$S(n)/a(n) \longrightarrow Y \qquad (n \to \infty). \tag{1.17}$$

We extend $a(n)$ into a function $a(t)$ on $\mathbb{R}^+$ by interpolation. By the classic theory of the central limit theorem, (1.4) holds also in this case.

Let the integer $p \geq 1$ be fixed and satisfy (1.5) and let $S_1(n), \dots, S_p(n)$ be $p$ independent copies of $S(n)$. The multi-parameter process $\overline{S}(n_1, \dots, n_p)$ given by

$$\overline{S}(n_1, \dots, n_p) = S_1(n_1) + \cdots + S_p(n_p), \qquad n_1, \dots, n_p = 0, 1, 2, \dots,$$

is called the additive random walk generated by the random walks $S_1(n), \dots, S_p(n)$. Its local time $l(n, x)$ is defined as

$$l(n, x) = \sum_{k_1, \dots, k_p = 0}^{n} \mathbb{1}_{\{S_1(k_1) + \cdots + S_p(k_p) = x\}}, \qquad x \in \mathbb{Z}^d,\ n = 1, 2, \dots.$$

Let $b_n$ be a deterministic positive sequence satisfying

$$b_n \longrightarrow \infty \quad \text{and} \quad b_n = o(n) \qquad (n \to \infty). \tag{1.18}$$

THEOREM 1.2. *Under (1.5), (1.17) and (1.18),*

$$\begin{aligned} &\lim_{n\to\infty} \frac{1}{b_n} \log \mathbb{P}\left\{ l(n, 0) \geq \lambda n^p a\left(\frac{n}{b_n}\right)^{-d} \right\} \\ &\quad = -(2\pi)^\alpha \frac{d}{\alpha}\left(1 - \frac{d}{\alpha p}\right)^{(\alpha p - d)/d} \left(\frac{\lambda}{\rho_1}\right)^{\alpha/d}, \end{aligned} \tag{1.19}$$

$$\begin{aligned} &\lim_{n\to\infty} \frac{1}{b_n} \log \mathbb{P}\left\{ \sum_{x \in \mathbb{Z}^d} l^2(n, x) \geq \lambda n^{2p} a\left(\frac{n}{b_n}\right)^{-d} \right\} \\ &\quad = -(2\pi)^\alpha \frac{d}{2\alpha}\left(1 - \frac{d}{2\alpha p}\right)^{(2\alpha p - d)/d} \left(\frac{\lambda}{\rho_2}\right)^{\alpha/d} \end{aligned} \tag{1.20}$$

*for every* $\lambda > 0$.

Theorems 1.1 and 1.2 apply to the following laws of the iterated logarithm.

THEOREM 1.3. *In the assumptions of Theorem 1.1,*

$$\begin{aligned} &\limsup_{t\to\infty} \frac{1}{t^p} a\left(\frac{t}{\log\log t}\right)^d L(t, 0) \\ &\quad = (2\pi)^{-d} \left(\frac{\alpha}{d}\right)^{d/\alpha} \left(1 - \frac{d}{\alpha p}\right)^{-(p - d/\alpha)} \rho_1 \qquad a.s., \end{aligned} \tag{1.21}$$

$$\limsup_{t\to\infty}\frac{1}{t^{2p}}a\biggl(\frac{t}{\log\log t}\biggr)^d I_t \tag{1.22}$$
$$= \biggl(\frac{2\alpha}{d}\biggr)^{d/\alpha}(2\pi)^{-d}\biggl(1-\frac{d}{2\alpha p}\biggr)^{-(2\alpha p-d)/\alpha}\rho_2 \qquad a.s.,$$

*where $I_t$ is given in* (1.15).

*In the assumptions of Theorem* 1.2,

$$\limsup_{n\to\infty}\frac{1}{n^{p}}a\biggl(\frac{n}{\log\log n}\biggr)^d l(n,0) \tag{1.23}$$
$$= (2\pi)^{-d}\biggl(\frac{\alpha}{d}\biggr)^{d/\alpha}\biggl(1-\frac{d}{\alpha p}\biggr)^{-(p-d/\alpha)}\rho_1 \qquad a.s.,$$

$$\limsup_{n\to\infty}\frac{1}{n^{2p}}a\biggl(\frac{n}{\log\log n}\biggr)^d \sum_{x\in\mathbb{Z}^d} l^2(n,x) \tag{1.24}$$
$$= (2\pi)^{-d}\biggl(\frac{2\alpha}{d}\biggr)^{d/\alpha}\biggl(1-\frac{d}{2\alpha p}\biggr)^{-(2\alpha p-d)/\alpha}\rho_2 \qquad a.s.$$

A special case covered by our theorems is when $X(t)$, $S(n)$ are square integrable and therefore are attracted, according to the classic central limit theorem, by the normal distributions with $a(t)=\sqrt{t}$. Let $\Gamma$ be the covariance matrix of $X(1)$ [or of $S(1)$]. By our assumption, $\Gamma$ is positive definite and

$$\Psi(\lambda)=\tfrac{1}{2}\lambda\cdot(\Gamma\lambda), \qquad \lambda\in\mathbb{R}^d.$$

Let $A$ be a $d\times d$ positive definite matrix such that $\Gamma=A^2$. Notice that, under the substitution $f(\lambda)=\sqrt[4]{\det\Gamma}\,g(A\lambda)$,

$$\begin{aligned}
&\int_{\mathbb{R}^d}\biggl[\int_{\mathbb{R}^d}\frac{f(\lambda+\gamma)f(\gamma)}{\sqrt{1+\Psi(\lambda+\gamma)}\sqrt{1+\Psi(\gamma)}}\,d\gamma\biggr]^p d\lambda\\
&\qquad=[\det\Gamma]^{p/2}\int_{\mathbb{R}^d}\biggl[\int_{\mathbb{R}^d}\frac{g(A(\lambda+\gamma))g(A\gamma)}{\sqrt{1+2^{-1}|A(\lambda+\gamma)|^2}\sqrt{1+2^{-1}|A\gamma|^2}}\,d\gamma\biggr]^p d\lambda\\
&\qquad=\frac{1}{\sqrt{\det\Gamma}}\int_{\mathbb{R}^d}\biggl[\int_{\mathbb{R}^d}\frac{g(\lambda+\gamma)g(\gamma)}{\sqrt{1+2^{-1}|\lambda+\gamma|^2}\sqrt{1+2^{-1}|\gamma|^2}}\,d\gamma\biggr]^p d\lambda,
\end{aligned}$$

where the second step follows from the linear transform $(\lambda,\gamma)\mapsto(A\lambda,A\gamma)$. Similarly,

$$\begin{aligned}
&\int_{\mathbb{R}^d}\biggl[\int_{\mathbb{R}^d}\frac{f(\lambda+\gamma)f(\gamma)}{\sqrt{1+\Psi(\lambda+\gamma)}\sqrt{1+\Psi(\gamma)}}\,d\gamma\biggr]^{2p} d\lambda\\
&\qquad=\frac{1}{\sqrt{\det\Gamma}}\int_{\mathbb{R}^d}\biggl[\int_{\mathbb{R}^d}\frac{g(\lambda+\gamma)g(\gamma)}{\sqrt{1+2^{-1}|\lambda+\gamma|^2}\sqrt{1+2^{-1}|\gamma|^2}}\,d\gamma\biggr]^{2p} d\lambda.
\end{aligned}$$

Hence, we have

$$\rho_1 = \frac{1}{\sqrt{\det\Gamma}}\bar\rho_1 \quad \text{and} \quad \rho_2 = \frac{1}{\sqrt{\det\Gamma}}\bar\rho_2,$$

where

$$\bar\rho_1 = \sup_{\|g\|_2=1}\int_{\mathbb{R}^d}\left[\int_{\mathbb{R}^d}\frac{g(\lambda+\gamma)g(\gamma)}{\sqrt{1+2^{-1}|\lambda+\gamma|^2}\sqrt{1+2^{-1}|\gamma|^2}}\,d\gamma\right]^p d\lambda,$$

$$\bar\rho_2 = \sup_{\|g\|_2=1}\int_{\mathbb{R}^d}\left[\int_{\mathbb{R}^d}\frac{g(\lambda+\gamma)g(\gamma)}{\sqrt{1+2^{-1}|\lambda+\gamma|^2}\sqrt{1+2^{-1}|\gamma|^2}}\,d\gamma\right]^{2p} d\lambda.$$

To avoid repeating statements, we only consider the additive random walk in the following corollary.

COROLLARY 1.4. *Assume that* $d < 2p$ *and assume that* $\mathbb{E}|S(1)|^2 < \infty$. *Then*

$$\begin{aligned}\lim_{n\to\infty}\frac{1}{b_n}\log\mathbb{P}\{l(n,0)\ge\lambda n^{(2p-d)/2}b_n^{d/2}\}\\ = -2d\pi^2\left(1-\frac{d}{2p}\right)^{(2p-d)/d}(\det\Gamma)^{1/d}\left(\frac{\lambda}{\bar\rho_1}\right)^{2/d},\end{aligned}\tag{1.25}$$

$$\begin{aligned}\lim_{n\to\infty}\frac{1}{b_n}\log\mathbb{P}\left\{\sum_{x\in\mathbb{Z}^d}l^2(n,x)\ge\lambda n^{(4p-d)/2}b_n^{d/2}\right\}\\ = -d\pi^2\left(1-\frac{d}{4p}\right)^{(4p-d)/d}(\det\Gamma)^{1/d}\left(\frac{\lambda}{\bar\rho_2}\right)^{2/d},\end{aligned}\tag{1.26}$$

$$\begin{aligned}\limsup_{n\to\infty} n^{-(2p-d)/2}(\log\log n)^{-d/2}l(n,0)\\ = (\sqrt{2d}\pi)^{-d}\left(1-\frac{d}{2p}\right)^{-(2p-d)/2}\frac{\bar\rho_1}{\sqrt{\det\Gamma}}\qquad a.s.,\end{aligned}\tag{1.27}$$

$$\begin{aligned}\limsup_{n\to\infty} n^{-(4p-d)/2}(\log\log n)^{-d/2}\sum_{x\in\mathbb{Z}^d}l^2(n,x)\\ = (\sqrt{d}\pi)^{-d}\left(1-\frac{d}{4p}\right)^{-(4p-d)/2}\frac{\bar\rho_2}{\sqrt{\det\Gamma}}\qquad a.s.\end{aligned}\tag{1.28}$$

REMARK. From Lemma 5.1 and (6.4) given later, one can see that if we replace $L(t,0)$ by $L(t,x)$ for a fixed $x$ in Theorems 1.1–1.3 and Corollary 1.4, then our results still hold.

The limit laws for the local times in the classic case $p=1$ have been extensively studied and it is impossible for us to list all works in the literature. We mention [3, 6, 14, 17, 20, 24, 28, 39, 41] for classic reference.

The study of the quadratic form of the local time is linked to the limit theorems for random walks in random sceneries. We refer the reader to [13, 29, 36, 42]. The study is also motivated by the needs from physics for investigating self-intersection of the random paths. We cite [1, 2, 11, 12] for some recent results on the large and moderate deviations for the self-intersection local times in the case $p=1$ and [19, 43, 44, 45, 46] for the physicists' view on the self-intersection local times. To see how our theorems connect to the problem of self-intersection, we take the additive random walks as example. We introduce the notation

$$\{(k_1,\dots,k_p),(l_1,\dots,l_p)\}$$

for a two-element set where each of the elements is a point in $(\mathbb{Z}^+)^p$. In particular, this set is viewed as being identical to the set

$$\{(l_1,\dots,l_p),(k_1,\dots,k_p)\}.$$

The object of interest is the random quantity

$$\begin{aligned}\Lambda_n = \#\{\{(k_1,\dots,k_p),(l_1,\dots,l_p)\}; 0\le k_1,\dots,k_p,l_1,\dots,l_p\le n \\ (k_1,\dots,k_p)\ne(l_1,\dots,l_p) \text{ and} \\ S_1(k_1)+\cdots+S_p(k_p) = S_1(l_1)+\cdots+S_p(l_p)\},\end{aligned}$$

which counts the self-intersection of the additive random walk $\overline{S}(n_1,\dots,n_p)$ during the "period" $[0,n]^p$. On the other hand, notice that

$$\begin{aligned}&\sum_{x\in\mathbb{Z}^d} l^2(n,x)\\ &\quad = \#\{(k_1,\dots,k_p,l_1,\dots,l_p)\in[0,n]^{2p};\\ &\qquad S_1(k_1)+\cdots+S_p(k_p) = S_1(l_1)+\cdots+S_p(l_p)\}\\ &\quad = 2\Lambda_n + n^p\end{aligned}$$

and that under (1.5),

$$n^p = o\left\{n^{2p} a\left(\frac{n}{b_n}\right)^{-d}\right\} \qquad (n\to\infty)$$

for every sequence $\{b_n\}$ satisfying (1.18). From (1.20) and (1.24), we have the following:

Corollary 1.5. *In the assumption of Theorem 1.2,*

$$\lim_{n\to\infty}\frac{1}{b_n}\log\mathbb{P}\left\{\Lambda_n\ge\lambda n^{2p}a\left(\frac{n}{b_n}\right)^{-d}\right\} \tag{1.29}$$
$$= -2^{-(\alpha+d)/d}(2\pi)^{\alpha}\frac{d}{\alpha}\left(1-\frac{d}{2\alpha p}\right)^{(2\alpha p-d)/d}\left(\frac{\lambda}{\rho_2}\right)^{\alpha/d},$$

$$\limsup_{n\to\infty}\frac{1}{n^{2p}}a\left(\frac{n}{\log\log n}\right)^{d}\Lambda_n \tag{1.30}$$
$$= 2\left(\frac{2\alpha}{d}\right)^{d/\alpha}(2\pi)^{-d}\left(1-\frac{d}{2\alpha p}\right)^{-(2\alpha p-d)/\alpha}\rho_2 \qquad a.s.$$

From the technical view point, the multi-parameter case is quite different from the single-parameter case. In the multi-parameter case, our ability of using the Markov property is severely limited. Due to lack of scaling properties, our approach is fundamentally different from the one used in [9, 10] in the stable case, where time exponentiation is essential to the solution.

To outline our approach, recall a general theorem (Theorem 4, [8]) of Gärtner–Ellis type: Let $\{Z_\varepsilon\}$ be a family of nonnegative random variables and let $p\ge 1$ be an integer. Assume that, for any $\theta>0$, the following limit exists:

$$\lim_{\varepsilon\to 0^+}\varepsilon\log\sum_{m=0}^{\infty}\frac{(\theta\varepsilon^{-1})^m}{m!}(\mathbb{E}Z_\varepsilon^m)^{1/p}=\Lambda(\theta).$$

Write

$$I(\lambda)=p\sup_{\theta>0}\{\lambda^{1/p}\theta-\Lambda(\theta)\}.$$

By an argument of duality, we have that, for any $\theta_0>0$,

$$\sup_{\lambda>0}\left\{\lambda^{1/p}\theta_0-\frac{1}{p}I(\lambda)\right\}=\Lambda(\theta_0).$$

Assume further that, for any $\lambda_0>0$, there is a $\theta_0>0$ such that $\lambda_0$ is the unique maximizer of the function

$$\varphi(\lambda)=\lambda^{1/p}\theta_0-\frac{1}{p}I(\lambda).$$

Then

$$\lim_{\varepsilon\to 0^+}\varepsilon\log\mathbb{P}\{Z_\varepsilon\ge\lambda\}=-I(\lambda),\qquad \lambda>0.$$

To prove Theorem 1.1 and (1.20) in Theorem 1.2, therefore, it is sufficient to show

$$\lim_{t\to\infty}\frac{1}{b_t}\log\sum_{m=0}^\infty\frac{\theta^m}{m!}\left(\frac{b_t}{t}\right)^m a\left(\frac{t}{b_t}\right)^{md/p}(\mathbb{E}L^m(t,0))^{1/p}$$
$$= (2\pi)^{-\alpha d/(\alpha p-d)}\rho_1^{\alpha/(\alpha p-d)}\theta^{\alpha p/(\alpha p-d)}, \tag{1.31}$$

$$\lim_{t\to\infty}\frac{1}{b_t}\log\sum_{m=0}^\infty\frac{\theta^m}{m!}\left(\frac{b_t}{t}\right)^m a\left(\frac{t}{b_t}\right)^{md/2p}\{\mathbb{E}I_t^{m/2}\}^{1/p}$$
$$= (2\pi)^{-\alpha d/(2\alpha p-d)}\rho_2^{\alpha/(2\alpha p-d)}\theta^{2\alpha p/(2\alpha p-d)}, \tag{1.32}$$

$$\lim_{n\to\infty}\frac{1}{b_n}\log\sum_{m=0}^\infty\frac{\theta^m}{m!}\left(\frac{b_n}{n}\right)^m a\left(\frac{n}{b_n}\right)^{md/2p}\left\{\mathbb{E}\left[\sum_{x\in\mathbb{Z}^d}l^2(n,x)\right]^{m/2}\right\}^{1/p}$$
$$= (2\pi)^{-\alpha d/(2\alpha p-d)}\rho_2^{\alpha/(2\alpha p-d)}\theta^{2\alpha p/(2\alpha p-d)}. \tag{1.33}$$

Unfortunately, the argument we shall use for (1.31) is broken in the case of additive random walks. The proof of (1.19) in Theorem 1.1 follows from a separate treatment. The basic tool we adapt to achieve these goals is Fourier transformation.

For technical reasons, the proof of our theorems is not given in the order that the theorems are stated. In Section 2 we prove the lower bounds for (1.31), (1.32) and (1.33). Our proof relies on the Feynman–Kac type minoration developed in [7]. In Section 3 we prove that the local times and the quadratic forms of the local times weakly converges to their stable counterpart if properly normalized. In addition to being interesting for its own sake, this result is needed in the proof of the upper bounds of (1.31), (1.32) and (1.33). In Section 4 we establish the upper bounds for (1.31), (1.32) and (1.33). The central pieces of this section are the establishment of some moment inequalities. By the time we end Section 4, the proof of Theorem 1.1 and (1.20) in Theorem 1.2 will be complete. In Section 5 we prove the laws of the iterated logarithm given in (1.21), (1.22) and (1.24) as the applications of Theorems 1.1 and 1.2. In Section 6 we prove (1.19) in Theorem 1.2 and (1.23) in Theorem 1.3. The basic idea is to approximate the local time of the additive random walks by the local time of a properly constructed additive Lévy process and the involved techniques include randomization, symmetrization and moment comparison. In the Appendix we prove an analytic lemma.

## 2. Lower bounds in (1.31), (1.32) and (1.33).

In this section we establish

$$\liminf_{t\to\infty}\frac{1}{b_t}\log\sum_{m=0}^\infty\frac{\theta^m}{m!}\left(\frac{b_t}{t}\right)^m a\left(\frac{t}{b_t}\right)^{md/p}(\mathbb{E}L^m(t,0))^{1/p}$$

$$\geq (2\pi)^{-\alpha d(p+1)/(\alpha p-d)}\rho_1^{\alpha/(\alpha p-d)}\theta^{\alpha p/(\alpha p-d)}, \tag{2.1}$$

$$\liminf_{t\to\infty}\frac{1}{b_t}\log\sum_{m=0}^{\infty}\frac{\theta^m}{m!}\left(\frac{b_t}{t}\right)^m a\left(\frac{t}{b_t}\right)^{md/2p}\{\mathbb{E}I_t^{m/2}\}^{1/p}$$
$$\geq (2\pi)^{-\alpha d/(2\alpha p-d)}\rho_2^{\alpha/(2\alpha p-d)}\theta^{2\alpha p/(2\alpha p-d)}, \tag{2.2}$$

$$\liminf_{n\to\infty}\frac{1}{b_n}\log\sum_{m=0}^{\infty}\frac{\theta^m}{m!}\left(\frac{b_n}{n}\right)^m a\left(\frac{n}{b_n}\right)^{md/2p}\left\{\mathbb{E}\left[\sum_{x\in\mathbb{Z}^d}l^2(n,x)\right]^{m/2}\right\}^{1/p}$$
$$\geq (2\pi)^{-\alpha d/(2\alpha p-d)}\rho_2^{\alpha/(2\alpha p-d)}\theta^{2\alpha p/(2\alpha p-d)}. \tag{2.3}$$

We start with a simple lemma.

LEMMA 2.1. *Let $X(t)$ be a symmetric Lévy process. For any $\lambda_1,\dots,\lambda_l\in\mathbb{R}^d$ and any $t_1,\dots,t_l\geq 0$,*

$$\mathbb{E}\exp\left\{i\sum_{k=1}^{l}\lambda_k\cdot X(t_k)\right\}>0.$$

PROOF. Without loss, generality, we may assume that $t_1\leq t_2\leq\cdots\leq t_l$:

$$\begin{aligned}\mathbb{E}\exp\left\{i\sum_{k=1}^{l}\lambda_k\cdot X(t_k)\right\}&=\mathbb{E}\exp\left\{i\sum_{k=1}^{l}\left(\sum_{j=k}^{l}\lambda_j\right)\cdot(X(t_k)-X(t_{k-1}))\right\}\\&=\prod_{k=1}^{l}\exp\left\{-(t_k-t_{k-1})\psi\left(\sum_{j=k}^{l}\lambda_j\right)\right\}>0,\end{aligned}$$

where we use the convention $t_0=0$. □

In the rest of the paper, the Fourier transformation will be frequently used. For any function $g(x)$ on $\mathbb{R}^d$, we use $\widehat{g}(\lambda)$ for its Fourier transform:

$$\widehat{g}(\lambda)=\int_{\mathbb{R}^d}g(x)e^{i\lambda\cdot x}\,dx,\qquad \lambda\in\mathbb{R}^d.$$

We refer the reader to the books by Edwards [21, 22] for the general information on Fourier analysis.

Write

$$\mathcal{F}_\Psi=\left\{g\in\mathcal{L}^2(\mathbb{R}^d);\int_{\mathbb{R}^d}g^2(x)\,dx=1\text{ and }\int_{\mathbb{R}^d}|\widehat{g}(\lambda)|^2\Psi(\lambda)\,d\lambda<\infty\right\}. \tag{2.4}$$

Finally, a function $f(x)$ on $\mathbb{R}^d$, $\mathbb{Z}^d$ or $[-\pi,\pi]^d$ is said to be symmetric if

$$f(-\lambda)=f(\lambda), \qquad \lambda\in\mathbb{R}^d,\ \mathbb{Z}^d \text{ or } [-\pi,\pi]^d.$$

LEMMA 2.2. *Let $f(\lambda)$ be a symmetric function on $\mathbb{R}^d$ such that*

$$\int_{\mathbb{R}^d}|f(\lambda)|\,d\lambda<\infty \tag{2.5}$$

*and write*

$$\bar f(x)=\int_{\mathbb{R}^d}f(\lambda)e^{i\lambda\cdot x}\,d\lambda, \qquad x\in\mathbb{R}^d.$$

*Let $\theta>0$ be fixed but arbitrary.*

(a) *In the assumptions of Theorem 1.1 and under (1.6),*

$$\begin{aligned}\liminf_{t\to\infty}\frac{1}{b_t}\log\sum_{m=0}^{\infty}\frac{\theta^m}{m!}\left(\frac{b_t}{t}\right)^m\int_{(\mathbb{R}^d)^m}d\lambda_1\cdots d\lambda_m\\ \times\left[\prod_{k=1}^m f(\lambda_k)\right]\left[\mathbb{E}\prod_{k=1}^m\int_0^t\exp\left\{ia\left(\frac{t}{b_t}\right)^{-1}\lambda_k\cdot X(s)\right\}ds\right]\\ \ge\sup_{g\in\mathcal{F}_\Psi}\left\{\theta\int_{\mathbb{R}^d}\bar f(x)g^2(x)\,dx-\int_{\mathbb{R}^d}|\widehat g(\lambda)|^2\Psi(\lambda)\,d\lambda\right\}.\end{aligned}\tag{2.6}$$

(b) *In the assumptions of Theorem 1.1 and under (1.7),*

$$\begin{aligned}\liminf_{t\to\infty}\frac{1}{b_t}\log\sum_{m=0}^{\infty}\frac{\theta^m}{m!}\left(\frac{b_t}{t}\right)^m\int_{(a(tb_t^{-1})[-\pi,\pi]^d)^m}d\lambda_1\cdots d\lambda_m\\ \times\left[\prod_{k=1}^m f(\lambda_k)\right]\left[\mathbb{E}\prod_{k=1}^m\int_0^t\exp\left\{ia\left(\frac{t}{b_t}\right)^{-1}\lambda_k\cdot X(s)\right\}ds\right]\\ \ge\sup_{g\in\mathcal{F}_\Psi}\left\{\theta\int_{\mathbb{R}^d}\bar f(x)g^2(x)\,dx-\int_{\mathbb{R}^d}|\widehat g(\lambda)|^2\Psi(\lambda)\,d\lambda\right\}.\end{aligned}\tag{2.7}$$

(c) *In the assumptions of Theorem 1.2,*

$$\begin{aligned}\liminf_{n\to\infty}\frac{1}{b_n}\log\sum_{m=0}^{\infty}\frac{\theta^m}{m!}\left(\frac{b_n}{n}\right)^m\int_{(a(nb_n^{-1})[-\pi,\pi]^d)^m}d\lambda_1\cdots d\lambda_m\\ \times\left(\prod_{k=1}^m f(\lambda_k)\right)\left[\mathbb{E}\prod_{k=1}^m\sum_{l=1}^n\exp\left\{ia\left(\frac{n}{b_n}\right)^{-1}\lambda_k\cdot S(l)\right\}\right]\\ \ge\sup_{g\in\mathcal{F}_\Psi}\left\{\theta\int_{\mathbb{R}^d}\bar f(x)g^2(x)\,dx-\int_{\mathbb{R}^d}|\widehat g(\lambda)|^2\Psi(\lambda)\,d\lambda\right\}.\end{aligned}\tag{2.8}$$

PROOF. Due to similarity, we only prove (2.6) under (1.6). For each integer $m \geq 1$,

$$\int_{(\mathbb{R}^d)^m} d\lambda_1 \cdots d\lambda_m \left[\prod_{k=1}^m f(\lambda_k)\right] \times \left[\mathbb{E} \prod_{k=1}^m \int_0^t \exp\left\{ia\left(\frac{t}{b_t}\right)^{-1} \lambda_k \cdot X(s)\right\} ds\right]$$
$$= \mathbb{E}\left\{\int_0^t \left[\int_{\mathbb{R}^d} f(\lambda) \exp\left\{ia\left(\frac{t}{b_t}\right)^{-1} \lambda \cdot X(s)\right\} d\lambda\right] ds\right\}^m$$
$$= \mathbb{E}\left[\int_0^t \bar{f}\left(a\left(\frac{t}{b_t}\right)^{-1} X(s)\right) ds\right]^m.$$

Hence,

$$\sum_{m=0}^\infty \frac{\theta^m}{m!}\left(\frac{b_t}{t}\right)^m \int_{(\mathbb{R}^d)^m} d\lambda_1 \cdots d\lambda_m \left(\prod_{k=1}^m f(\lambda_k)\right) \times \left[\mathbb{E} \prod_{k=1}^m \int_0^t \exp\left\{ia\left(\frac{t}{b_t}\right)^{-1} \lambda_k \cdot X(s)\right\} ds\right]$$
$$= \mathbb{E} \exp\left\{\theta \frac{b_t}{t} \int_0^t \bar{f}\left(a\left(\frac{t}{b_t}\right)^{-1} X(s)\right) ds\right\}.$$

Under (2.5), $\bar{f}(x)$ is a real, bounded and continuous function on $\mathbb{R}^d$. By the same argument as Theorem 4.1 of [7],

$$\liminf_{t\to\infty} \frac{1}{b_t} \log \mathbb{E} \exp\left\{\theta \frac{b_t}{t} \int_0^t \bar{f}\left(a\left(\frac{t}{b_t}\right)^{-1} X(s)\right) ds\right\}$$
$$\geq \sup_{g \in \mathcal{F}_\Psi} \left\{\theta \int_{\mathbb{R}^d} \bar{f}(x) g^2(x)\, dx - \int_{\mathbb{R}^d} |\widehat{g}(\lambda)|^2 \Psi(\lambda)\, d\lambda\right\}.$$

This leads to (2.6). □

We are back to the proof of (2.1), (2.2) and (2.3). Due to similarity, we only prove (2.1) and (2.2) under the condition (1.6). For any nonnegative and symmetric function $f(\lambda)$ satisfying (2.5), write

$$\rho(f) = \sup_{\|g\|_2=1} \iint_{\mathbb{R}^d \times \mathbb{R}^d} f(\lambda - \gamma) \frac{g(\lambda) g(\gamma)}{\sqrt{1+\Psi(\lambda)}\sqrt{1+\Psi(\gamma)}}\, d\lambda\, d\gamma$$
$$= \sup_{\|g\|_2=1} \int_{\mathbb{R}^d} f(\lambda) \left[\int_{\mathbb{R}^d} \frac{g(\lambda+\gamma) g(\gamma)}{\sqrt{1+\Psi(\lambda+\gamma)}\sqrt{1+\Psi(\gamma)}}\, d\gamma\right] d\lambda.$$

We now prove (2.1). By inverse Fourier transformation,

$$
\begin{aligned}
L(t,0) &= \frac{1}{(2\pi)^d}\int_{\mathbb{R}^d}\left[\int_{\mathbb{R}^d} e^{-i\lambda\cdot x}L(t,x)\,dx\right]d\lambda\\
&= \frac{1}{(2\pi)^d}\int_{\mathbb{R}^d}\left[\int_0^t\cdots\int_0^t \exp\left\{-i\lambda\cdot(X_1(s_1)+\cdots\right.\right.\\
&\qquad\qquad\left.\left.+X_p(s_p))\right\}ds_1\cdots ds_p\right]d\lambda\\
&= \frac{1}{(2\pi)^d}\int_{\mathbb{R}^d} d\lambda\prod_{j=1}^p\int_0^t e^{-i\lambda\cdot X_j(s)}\,ds.
\end{aligned}
\tag{2.9}
$$

For any integer $m\ge 1$,

$$
\mathbb{E}L^m(t,0) = (2\pi)^{-md}\int_{(\mathbb{R}^d)^m} d\lambda_1\cdots d\lambda_m\left[\mathbb{E}\prod_{k=1}^m\int_0^t e^{i\lambda_k\cdot X(s)}\,ds\right]^p.
\tag{2.10}
$$

Let $q>1$ be the conjugate number of $p$ defined by the relation $p^{-1}+q^{-1}=1$. We now let $f(\lambda)$ be a symmetric nonnegative function on $\mathbb{R}^d$ satisfying (2.5) and $\|f\|_q=1$. By rescaling,

$$
\begin{aligned}
&(\mathbb{E}L^m(t,0))^{1/p}\\
&\quad= (2\pi)^{-md/p}a\left(\frac{t}{b_t}\right)^{-md/p}\\
&\qquad\times\left\{\int_{(\mathbb{R}^d)^m} d\lambda_1\cdots d\lambda_m\right.\\
&\qquad\qquad\left.\times\left[\mathbb{E}\prod_{k=1}^m\int_0^t\exp\left\{ia\left(\frac{t}{b_t}\right)^{-1}\lambda_k\cdot X(s)\right\}ds\right]^p\right\}^{1/p}\\
&\quad\ge (2\pi)^{-md/p}a\left(\frac{t}{b_t}\right)^{-md/p}\int_{(\mathbb{R}^d)^m} d\lambda_1\cdots d\lambda_m\\
&\qquad\times\left[\prod_{k=1}^m f(\lambda_k)\right]\left[\mathbb{E}\prod_{k=1}^m\int_0^t\exp\left\{ia\left(\frac{t}{b_t}\right)^{-1}\lambda_k\cdot X(s)\right\}ds\right],
\end{aligned}
$$

where the inequality follows from the Hölder inequality and from the crucial fact that

$$
\mathbb{E}\prod_{k=1}^m\int_0^t\exp\left\{ia\left(\frac{t}{b_t}\right)^{-1}\lambda_k\cdot X(s)\right\}ds\ge 0,
\tag{2.11}
$$

which is supported by Lemma 2.1.

Taking $\theta = \rho(f)^{-1}(2\pi)^{d/p}$ in Lemma 2.2 gives

$$\begin{aligned}
&\liminf_{t\to\infty} \frac{1}{b_t} \log \sum_{m=0}^{\infty} \frac{1}{m!} (\rho(f)^{-1}(2\pi)^{d/p})^m \\
&\qquad \times \left(\frac{b_t}{t}\right)^m a\left(\frac{t}{b_t}\right)^{md/p} (\mathbb{E}L^m(t,0))^{1/p} \qquad (2.12)\\
&\qquad \ge \sup_{g\in\mathcal{F}_\Psi} \left\{ \frac{1}{\rho(f)} \int_{\mathbb{R}^d} \bar{f}(x) g^2(x)\,dx - \int_{\mathbb{R}^d} |\widehat{g}(\lambda)|^2 \Psi(\lambda)\,d\lambda \right\} = 1,
\end{aligned}$$

where the last step follows from Lemma A.1 in the Appendix.

Let $\theta > 0$ be given as in (2.6) and let $\varepsilon > 0$ be arbitrarily small. By regularity of $a(t)$ given in (1.4), for large $t$, we have

$$\begin{aligned}
&\sum_{m=0}^{\infty} \frac{\theta^m}{m!} \left(\frac{b_t}{t}\right)^m a\left(\frac{t}{b_t}\right)^{md/p} (\mathbb{E}L^m(t,0))^{1/p} \\
&\quad = \sum_{m=0}^{\infty} \frac{1}{m!} (\rho(f)^{-1}(2\pi)^{d/p})^m (\theta\rho(f)(2\pi)^{-d/p})^m \\
&\qquad \times \left(\frac{b_t}{t}\right)^m a\left(\frac{t}{b_t}\right)^{md/p} (\mathbb{E}L^m(t,0))^{1/p} \\
&\quad \ge \sum_{m=0}^{\infty} \frac{1}{m!} (\rho(f)^{-1}(2\pi)^{d/p})^m \left(\frac{\lambda b_t}{t}\right)^m a\left(\frac{t}{\lambda b_t}\right)^{md/p} (\mathbb{E}L^m(t,0))^{1/p},
\end{aligned}$$

where

$$\lambda = (1-\varepsilon)(\theta\rho(f)(2\pi)^{-d/p})^{\alpha p/(\alpha p - d)}.$$

Replacing $b_t$ by $\lambda b_t$ in (2.12),

$$\begin{aligned}
&\liminf_{t\to\infty} \frac{1}{b_t} \log \sum_{m=0}^{\infty} \frac{\theta^m}{m!} \left(\frac{b_t}{t}\right)^m a\left(\frac{t}{b_t}\right)^{md/p} (\mathbb{E}L^m(t,0))^{1/p} \\
&\qquad \ge (1-\varepsilon)(\theta\rho(f)(2\pi)^{-d/p})^{\alpha p/(\alpha p - d)}.
\end{aligned}$$

Letting $\varepsilon \to 0^+$,

$$\begin{aligned}
&\liminf_{t\to\infty} \frac{1}{b_t} \log \sum_{m=0}^{\infty} \frac{\theta^m}{m!} \left(\frac{b_t}{t}\right)^m a\left(\frac{t}{b_t}\right)^{md/p} (\mathbb{E}L^m(t,0))^{1/p} \qquad (2.13)\\
&\qquad \ge (2\pi)^{-\alpha d/(\alpha p - d)} \rho(f)^{\alpha p/(\alpha p - d)} \theta^{\alpha p/(\alpha p - d)}.
\end{aligned}$$

Notice that, for any $g \in \mathcal{L}^2(\mathbb{R}^d)$ with $g \ge 0$ and $\|g\|_2 = 1$, the function

$$H_g(\lambda) \equiv \int_{\mathbb{R}^d} \frac{g(\lambda+\gamma)g(\gamma)}{\sqrt{1+\Psi(\lambda+\gamma)}\sqrt{1+\Psi(\gamma)}}\,d\gamma$$

is nonnegative and symmetric on $\mathbb{R}^d$. Hence,

$$\begin{aligned}\sup_f \rho(f) &\ge \sup_f \sup_{\substack{\|g\|_2=1\\ g\ge 0}} \int_{\mathbb{R}^d} f(\lambda) H_g(\lambda)\, d\lambda\\ &= \sup_{\substack{\|g\|_2=1\\ g\ge 0}} \sup_f \int_{\mathbb{R}^d} f(\lambda) H_g(\lambda)\, d\lambda\\ &= \sup_{\substack{\|g\|_2=1\\ g\ge 0}} \left(\int_{\mathbb{R}^d} |H_g(\lambda)|^p\, d\lambda\right)^{1/p} = \rho_1^{1/p},\end{aligned}$$

where "$\sup_f$" is taken over all symmetric $f$ satisfying (2.5) with $\|f\|_q = 1$.

Summarizing what we have proved since (2.13), we have (2.1).

We now come to the proof of (2.2) under (1.6). By Parseval's identity, for any $t > 0$,

$$\begin{aligned}&\int_{\mathbb{R}^d} L^2(t,x)\, dx\\ &\qquad = \frac{1}{(2\pi)^d} \int_{\mathbb{R}^d} \left| \int_{\mathbb{R}^d} L(t,x) e^{i\lambda\cdot x}\, dx\right|^2 d\lambda\\ &\qquad = \frac{1}{(2\pi)^d} \int_{\mathbb{R}^d} \left| \int_{[0,t]^p} \exp\{i\lambda\cdot(X_1(s_i) + \cdots\right.\\ &\qquad\qquad\qquad\qquad \left. + X_p(s_p))\}\, ds_1\cdots ds_p\right|^2 d\lambda \qquad (2.14)\\ &\qquad = \frac{1}{(2\pi)^d} \int_{\mathbb{R}^d} \left| \prod_{j=1}^p \int_0^t e^{i\lambda\cdot X_j(s)}\, ds\right|^2 d\lambda.\end{aligned}$$

Let $f(\lambda) \ge 0$ be symmetric on $\mathbb{R}^d$ such that $\|f\|_2 = 1$. By rescaling and the Cauchy–Schwarz inequality,

$$\begin{aligned}&\left[\int_{\mathbb{R}^d} L^2(t,x)\, dx\right]^{1/2}\\ &\qquad \ge (2\pi)^{d/2} a \left(\frac{t}{b_t}\right)^{-d/2}\\ &\qquad\quad \times \left[\int_{\mathbb{R}^d} \left| \prod_{j=1}^p \int_0^t \exp\left\{ia\left(\frac{t}{b_t}\right)^{-1} \lambda\cdot X_j(s)\right\} ds\right|^2 d\lambda\right]^{1/2}\\ &\qquad \ge (2\pi)^{d/2} a \left(\frac{t}{b_t}\right)^{-d/2}\end{aligned}$$

$$\times \int_{\mathbb{R}^d} f(\lambda) \left| \prod_{j=1}^p \int_0^t \exp\left\{ ia\left(\frac{t}{b_t}\right)^{-1} \lambda \cdot X_j(s) \right\} ds \right| d\lambda.$$

For any $m \geq 1$,

$$\begin{aligned}
&\mathbb{E}\left[ \int_{\mathbb{R}^d} d\lambda\, f(\lambda) \left| \prod_{j=1}^p \int_0^t \exp\left\{ ia\left(\frac{t}{b_t}\right)^{-1} \lambda \cdot X_j(s) \right\} ds \right| \right]^m \\
&\quad = \int_{(\mathbb{R}^d)^m} d\lambda_1 \cdots d\lambda_m \left( \prod_{k=1}^m f(\lambda_k) \right) \\
&\qquad \times \left[ \mathbb{E} \left| \prod_{k=1}^m \int_0^t \exp\left\{ ia\left(\frac{t}{b_t}\right)^{-1} \lambda_k \cdot X(s) \right\} ds \right| \right]^p \\
&\quad \geq \int_{(\mathbb{R}^d)^m} d\lambda_1 \cdots d\lambda_m \left( \prod_{k=1}^m f(\lambda_k) \right) \\
&\qquad \times \left| \mathbb{E} \prod_{k=1}^m \int_0^t \exp\left\{ ia\left(\frac{t}{b_t}\right)^{-1} \lambda_k \cdot X(s) \right\} ds \right|^p.
\end{aligned}$$

Let $q > 1$ be the conjugate number of $p$ and let $h$ be a symmetric function on $\mathbb{R}^d$ satisfying

$$\int_{\mathbb{R}^d} |h(\lambda)| f(\lambda)\, d\lambda < \infty \quad \text{and} \quad \int_{\mathbb{R}^d} |h(\lambda)|^q f(\lambda)\, d\lambda = 1. \tag{2.15}$$

We have

$$\begin{aligned}
&\left\{ \mathbb{E}\left[ \int_{\mathbb{R}^d} L^2(t,x)\, dx \right]^{m/2} \right\}^{1/p} \\
&\quad \geq (2\pi)^{-md/(2p)} a\left(\frac{t}{b_t}\right)^{-md/(2p)} \\
&\qquad \times \left\{ \int_{(\mathbb{R}^d)^m} d\lambda_1 \cdots d\lambda_m \left( \prod_{k=1}^m f(\lambda_k) \right) \right. \\
&\qquad\quad \left. \times \left| \mathbb{E} \prod_{k=1}^m \int_0^t \exp\left\{ ia\left(\frac{t}{b_t}\right)^{-1} \lambda_k \cdot X(s) \right\} ds \right|^p \right\}^{1/p} \\
&\quad \geq (2\pi)^{-md/(2p)} a\left(\frac{t}{b_t}\right)^{-md/(2p)} \int_{(\mathbb{R}^d)^m} d\lambda_1 \cdots d\lambda_m \left( \prod_{k=1}^m (fh)(\lambda_k) \right) \\
&\qquad \times \left[ \mathbb{E} \prod_{k=1}^m \int_0^t \exp\left\{ ia\left(\frac{t}{b_t}\right)^{-1} \lambda_k \cdot X(s) \right\} ds \right].
\end{aligned}$$

Taking $\theta = (2\pi)^{d/2p}\rho(fh)^{-1}$ and replacing $f$ by $fh$ in (2.6) of Lemma 2.2,

$$\liminf_{t\to\infty}\frac{1}{b_t}\log\sum_{m=0}^{\infty}\frac{1}{m!}\left(\frac{(2\pi)^{d/2p}}{\rho(fh)}\right)^m \times\left(\frac{b_t}{t}\right)^m a\left(\frac{t}{b_t}\right)^{md/2p}\left\{\mathbb{E}\left[\int_{\mathbb{R}^d}L^2(t,x)\,dx\right]^{m/2}\right\}^{1/p} \geq \sup_{g\in\mathcal{F}_\Psi}\left\{\frac{1}{\rho(fh)}\int_{\mathbb{R}^d}\overline{fh}(x)g^2(x)\,dx-\int_{\mathbb{R}^d}|\widehat{g}(\lambda)|^2\Psi(\lambda)\,d\lambda\right\}=1,$$

where the last step follows from Lemma A.1 given in the Appendix.

As for the general $\theta>0$, similar to the argument used in the proof of (2.1),

$$\liminf_{t\to\infty}\frac{1}{b_t}\log\sum_{m=0}^{\infty}\frac{\theta^m}{m!}\left(\frac{b_t}{t}\right)^m a\left(\frac{t}{b_t}\right)^{md/2p} \times\left\{\mathbb{E}\left[\int_{\mathbb{R}^d}L^2(t,x)\,dx\right]^{m/2}\right\}^{1/p} \geq (2\pi)^{-\alpha d/(2\alpha p-d)}\rho(fh)^{2\alpha p/(2\alpha p-d)}\theta^{2\alpha p/(2\alpha p-d)}. \tag{2.16}$$

Notice that

$$\sup_h\rho(fh)=\sup_{\|g\|_2=1}\left\{\int_{\mathbb{R}^d}d\lambda\, f(\lambda)\left[\int_{\mathbb{R}^d}\frac{g(\lambda+\gamma)g(\gamma)}{\sqrt{1+\Psi(\lambda+\gamma)}\sqrt{1+\Psi(\gamma)}}\,d\gamma\right]^p\right\}^{1/p},$$

where the supremum is taken over all symmetric functions $h$ satisfying (2.15). Taking supremum on both sides over all symmetric, nonnegative functions $f$ with $\|f\|_2=1$, we have

$$\sup_{f,h}\rho(fh)=\rho_2^{1/2p}.$$

Finally, taking supremum over $h$ and $f$ on the right-hand side of (2.16) proves (2.2).

Remark. A careful reader may notice that, in the context of the additive random walks, the statement corresponding to (2.1) in the Lévy case is missing. The reason is the absence of the property like (2.11) in the case of the random walks. This also creates a problem in the proof of the upper bounds.

**3. Laws of weak convergence.** In our assumptions, the Lévy processes and random walks are attracted by stable processes. Naturally, we expect that this relation passes to the local times. Recall that $Y$ is a stable random variable given in (1.2) and (1.16), $Y(t)$ is a stable process in $\mathbb{R}^d$ such that $Y(1) \stackrel{d}{=} Y$, and $L_Y(t,x)$ is the local time of the additive stable process generated by $Y(t)$.

THEOREM 3.1. *In the assumptions of Theorem 1.1,*

$$\frac{a(t)^d}{t^p} L(t,0) \xrightarrow{d} L_Y(1,0) \qquad (t \to \infty), \tag{3.1}$$

$$\frac{a(t)^d}{t^{2p}} I_t \xrightarrow{d} \int_{\mathbb{R}^d} L_Y^2(1,x)\,dx \qquad (t \to \infty), \tag{3.2}$$

*where $I_t$ is given in (1.15).*

*In the assumptions of Theorem 1.2,*

$$\frac{a(n)^d}{n^p} l(n,0) \xrightarrow{d} L_Y(1,0) \qquad (n \to \infty), \tag{3.3}$$

$$\frac{a(n)^d}{n^{2p}} \sum_{x\in\mathbb{Z}^d} l^2(n,x) \xrightarrow{d} \int_{\mathbb{R}^d} L_Y^2(1,x)\,dx \qquad (n \to \infty). \tag{3.4}$$

PROOF. Due to similarity, we only prove (3.1) and (3.2) under the condition (1.6). We first rescale in (2.9) and (2.14):

$$L(t,0) = \frac{1}{(2\pi)^d} a(t)^{-d} \int_{\mathbb{R}^d} d\lambda \left[\prod_{j=1}^p \int_0^t \exp\left\{i\lambda \cdot \frac{X_j(s)}{a(t)}\right\} ds\right],$$

$$\int_{\mathbb{R}^d} L^2(t,x)\,dx = \frac{1}{(2\pi)^d} a(t)^{-d} \int_{\mathbb{R}^d} d\lambda \prod_{j=1}^p \left|\int_0^t \exp\left\{i\lambda \cdot \frac{X_j(s)}{a(t)}\right\} ds\right|^2.$$

Let $D\{[0,1],(\mathbb{R}^d)^p\}$ be the space of the $(\mathbb{R}^d)^p$-valued functions on $[0,1]$ which are right continuous and have left limits on $[0,1]$. Under the uniform convergence topology, $D\{[0,1],(\mathbb{R}^d)^p\}$ is a Banach space. For any fixed $M > 0$, the functionals

$$\mathcal{F}(x_1,\dots,x_p) = \int_{[-M,M]^d} d\lambda \left[\prod_{j=1}^p \int_0^1 e^{i\lambda\cdot x_j(s)}\,ds\right],$$

$$\mathcal{G}(x_1,\dots,x_p) = \int_{[-M,M]^d} d\lambda \prod_{j=1}^p \left|\int_0^1 e^{i\lambda\cdot x_j(s)}\,ds\right|^2$$

are continuous on $D\{[0,1],(\mathbb{R}^d)^p\}$. By the invariance principle,

$$\frac{1}{t^p} \int_{[-M,M]^d} d\lambda \left[\prod_{j=1}^p \int_0^t \exp\left\{i\lambda \cdot \frac{X_j(s)}{a(t)}\right\} ds\right]$$

$$\xrightarrow{d} \int_{[-M,M]^d} d\lambda \left[\prod_{j=1}^p \int_0^1 e^{i\lambda\cdot Y_j(s)}\,ds\right],$$

$$\frac{1}{t^{2p}}\int_{[-M,M]^d} d\lambda \prod_{j=1}^p \left|\int_0^t \exp\left\{i\lambda\cdot\frac{X_j(s)}{a(t)}\right\}ds\right|^2$$

$$\xrightarrow{d} \int_{[-M,M]^d} d\lambda \prod_{j=1}^p \left|\int_0^1 e^{i\lambda\cdot Y_j(s)}\,ds\right|^2,$$

as $t\to\infty$.

To prove (3.1) and (3.3), therefore, we need only to show that

$$\lim_{M\to\infty}\limsup_{t\to\infty}\frac{1}{t^{2p}}\mathbb{E}\left|\int_{\mathbb{R}^d\setminus[-M,M]^d} d\lambda\left[\prod_{j=1}^p\int_0^t \exp\left\{i\lambda\cdot\frac{X_j(s)}{a(t)}\right\}ds\right]\right|^2 = 0 \tag{3.5}$$

and that

$$\lim_{M\to\infty}\limsup_{t\to\infty}\frac{1}{t^{2p}}\int_{\mathbb{R}^d\setminus[-M,M]^d} d\lambda\,\mathbb{E}\prod_{j=1}^p\left|\int_0^t \exp\left\{i\lambda\cdot\frac{X_j(s)}{a(t)}\right\}ds\right|^2 = 0. \tag{3.6}$$

Notice that

$$\begin{aligned}
&\mathbb{E}\left|\int_{\mathbb{R}^d\setminus[-M,M]^d} d\lambda\left[\prod_{j=1}^p\int_0^t \exp\left\{i\lambda\cdot\frac{X_j(s)}{a(t)}\right\}ds\right]\right|^2\\
&\quad= \int_{(\mathbb{R}^d\setminus[-M,M]^d)^2} d\lambda\,d\gamma\left[\int_0^t\int_0^t \mathbb{E}\exp\left\{i\frac{\lambda\cdot X(s_1)-\gamma\cdot X(s_2)}{a(t)}\right\}ds_1\,ds_2\right]^p\\
&\quad\le 2^p\int_{(\mathbb{R}^d\setminus[-M,M]^d)^2} d\lambda\,d\gamma\left[\int_0^t\int_{s_1}^t \mathbb{E}\exp\left\{i\frac{\lambda\cdot X(s_1)-\gamma\cdot X(s_2)}{a(t)}\right\}ds_2\,ds_1\right]^p\\
&\quad= 2^p\int_{(\mathbb{R}^d\setminus[-M,M]^d)^2} d\lambda\,d\gamma\\
&\qquad\times\left[\int_0^t\int_{s_1}^t \exp\left\{-s_1\psi\left(\frac{\lambda-\gamma}{a(t)}\right)-(s_2-s_1)\psi\left(\frac{\gamma}{a(t)}\right)\right\}ds_2\,ds_1\right]^p\\
&\quad\le 2^p\int_{(\mathbb{R}^d\setminus[-M,M]^d)^2} d\lambda\,d\gamma\left[\int_0^t \exp\left\{-s\psi\left(\frac{\lambda-\gamma}{a(t)}\right)\right\}ds\right]^p\\
&\qquad\times\left[\int_0^t \exp\left\{-s\psi\left(\frac{\gamma}{a(t)}\right)\right\}ds\right]^p\\
&\quad\le 2^p\left\{\int_{\mathbb{R}^d} d\lambda\left[\int_0^t \exp\left\{-s\psi\left(\frac{\lambda}{a(t)}\right)\right\}ds\right]^p\right\}\\
&\qquad\times\left\{\int_{\mathbb{R}^d\setminus[-M,M]^d} d\gamma\left[\int_0^t \exp\left\{-s\psi\left(\frac{\gamma}{a(t)}\right)\right\}ds\right]^p\right\},
\end{aligned}$$

where the last step partially follows from the substitution $\lambda - \gamma \mapsto \lambda$ and $\gamma \mapsto \gamma$.

Hence, we need only to show that

$$\limsup_{t\to\infty} \frac{1}{t^p} \int_{\mathbb{R}^d} d\lambda \left[\int_0^t \exp\left\{-s\psi\left(\frac{\lambda}{a(t)}\right)\right\} ds\right]^p < \infty \tag{3.7}$$

and that

$$\lim_{M\to\infty} \limsup_{t\to\infty} \frac{1}{t^p} \int_{\mathbb{R}^d\setminus[-M,M]^d} d\gamma \left[\int_0^t \exp\left\{-s\psi\left(\frac{\gamma}{a(t)}\right)\right\} ds\right]^p = 0. \tag{3.8}$$

Indeed,

$$\begin{aligned}
&\frac{1}{t^p} \int_{\mathbb{R}^d\setminus[-M,M]^d} d\gamma \left[\int_0^t \exp\left\{-s\psi\left(\frac{\gamma}{a(t)}\right)\right\} ds\right]^p \\
&\qquad = \int_{\mathbb{R}^d\setminus[-M,M]^d} \left[t\psi\left(\frac{\gamma}{a(t)}\right)\right]^{-p} \left[1 - \exp\left\{-t\psi\left(\frac{\gamma}{a(t)}\right)\right\}\right]^p d\gamma.
\end{aligned}$$

By (1.5), there is a $\varepsilon > 0$ such that $(\alpha - \varepsilon)p > d$. By (1.6), there is a $N > 0$ such that $\psi(\lambda) \geq 2$ for all $\lambda \notin [-N, N]^d$. In view of (1.1), by an estimate similar to the one used for (2.i) in [38], there is a constant $\delta > 0$ such that

$$t\psi\left(\frac{\lambda}{a(t)}\right) \geq \delta|\lambda|^{\alpha-\varepsilon}, \qquad \lambda \in [-a(t)N, a(t)N]^d, \tag{3.9}$$

for large $t$. Thus,

$$\begin{aligned}
&\int_{\mathbb{R}^d\setminus[-M,M]^d} \left[t\psi\left(\frac{\gamma}{a(t)}\right)\right]^{-p} \left[1 - \exp\left\{-t\psi\left(\frac{\gamma}{a(t)}\right)\right\}\right]^p d\gamma \\
&\qquad \leq \delta^{-p} \int_{\mathbb{R}^d\setminus[-M,M]^d} |\lambda|^{-p(\alpha-\varepsilon)} d\lambda + \int_{\mathbb{R}^d\setminus[-a(t)N,a(t)N]^d} \left[t\psi\left(\frac{\gamma}{a(t)}\right)\right]^{-p} d\gamma.
\end{aligned}$$

By (1.4), (1.5) and (1.6),

$$\begin{aligned}
&\int_{\mathbb{R}^d\setminus[-a(t)N,a(t)N]^d} \left[t\psi\left(\frac{\gamma}{a(t)}\right)\right]^{-p} d\gamma \\
&\qquad = t^{-p} a(t)^d \int_{\mathbb{R}^d\setminus[-N,N]^d} [\psi(\lambda)]^{-p} d\lambda \longrightarrow 0
\end{aligned}$$

as $t \to \infty$. Hence, (3.8) holds. The proof of (3.7) is similar.

We now come to the proof of (3.6). Notice that

$$\begin{aligned}
&\frac{1}{t^{2p}} \int_{\mathbb{R}^d\setminus[-M,M]^d} d\lambda\, \mathbb{E} \prod_{j=1}^p \left|\int_0^t \exp\left\{i\lambda \cdot \frac{X_j(s)}{a(t)}\right\} ds\right|^2 \\
&\qquad = \frac{1}{t^{2p}} \int_{\mathbb{R}^d\setminus[-M,M]^d} d\lambda \left[\mathbb{E}\left|\int_0^t \exp\left\{i\lambda \cdot \frac{X(s)}{a(t)}\right\} ds\right|^2\right]^p
\end{aligned}$$

$$= \frac{1}{t^{2p}} \int_{\mathbb{R}^d \setminus [-M,M]^d} d\lambda \Big[ 2\mathbb{E} \int_0^t \int_{s_1}^t \exp\Big\{ i\lambda \cdot \frac{X(s_2) - X(s_1)}{a(t)} \Big\} ds_2\, ds_1 \Big]^p$$

and that

$$\begin{aligned}
&\mathbb{E} \int_0^t \int_{s_1}^t \exp\Big\{ i\lambda \cdot \frac{X(s_2) - X(s_1)}{a(t)} \Big\} ds_2\, ds_1 \\
&\quad = \int_0^t \int_{s_1}^t \exp\Big\{ -(s_2 - s_1)\psi\Big(\frac{\lambda}{a(t)}\Big) \Big\} ds_2\, ds_1 \\
&\quad = t\psi\Big(\frac{\lambda}{a(t)}\Big)^{-1} - \psi\Big(\frac{\lambda}{a(t)}\Big)^{-2} \Big[ 1 - \exp\Big\{ -t\psi\Big(\frac{\lambda}{a(t)}\Big) \Big\} \Big] \\
&\quad \le t\psi\Big(\frac{\lambda}{a(t)}\Big)^{-1}.
\end{aligned}$$

Hence,

$$\begin{aligned}
&\frac{1}{t^{2p}} \int_{\mathbb{R}^d \setminus [-M,M]^d} d\lambda\, \mathbb{E} \prod_{j=1}^p \Big| \int_0^t \exp\Big\{ i\lambda \cdot \frac{X_j(s)}{a(t)} \Big\} ds \Big|^2 \\
&\qquad \le 2^p \int_{\mathbb{R}^d \setminus [-M,M]^d} \Big[ t\psi\Big(\frac{\lambda}{a(t)}\Big) \Big]^{-p} d\lambda. 
\end{aligned} \tag{3.10}$$

By (3.9), (3.6) holds. □

**4. Upper bounds in (1.31), (1.32) and (1.33).** In this section we establish

$$\begin{aligned}
&\limsup_{t\to\infty} \frac{1}{b_t} \log \sum_{m=0}^\infty \frac{\theta^m}{m!} \Big(\frac{b_t}{t}\Big)^m a\Big(\frac{t}{b_t}\Big)^{md/p} (\mathbb{E} L^m(t,0))^{1/p} \\
&\qquad \le (2\pi)^{-\alpha d(p+1)/(\alpha p - d)} \rho_1^{\alpha/(\alpha p - d)} \theta^{\alpha p/(\alpha p - d)},
\end{aligned} \tag{4.1}$$

$$\begin{aligned}
&\limsup_{t\to\infty} \frac{1}{b_t} \log \sum_{m=0}^\infty \frac{\theta^m}{m!} \Big(\frac{b_t}{t}\Big)^m a\Big(\frac{t}{b_t}\Big)^{md/2p} \{\mathbb{E} I_t^{m/2}\}^{1/p} \\
&\qquad \le (2\pi)^{-\alpha d/(2\alpha p - d)} \rho_2^{\alpha/(2\alpha p - d)} \theta^{2\alpha p/(2\alpha p - d)},
\end{aligned} \tag{4.2}$$

$$\begin{aligned}
&\limsup_{n\to\infty} \frac{1}{b_n} \log \sum_{m=0}^\infty \frac{\theta^m}{m!} \Big(\frac{b_n}{n}\Big)^m a\Big(\frac{n}{b_n}\Big)^{md/2p} \Big\{ \mathbb{E} \Big[ \sum_{x\in\mathbb{Z}^d} l^2(n,x) \Big]^{m/2} \Big\}^{1/p} \\
&\qquad \le (2\pi)^{-\alpha d/(2\alpha p - d)} \rho_2^{\alpha/(2\alpha p - d)} \theta^{2\alpha p/(2\alpha p - d)}.
\end{aligned} \tag{4.3}$$

We first concentrate on the proof of (4.1) and (4.2) and then work (4.3) out later.

LEMMA 4.1. *Under both* (1.6) *and* (1.7), *for any* $t \geq 0$ *and any integer* $m \geq 1$,

$$\mathbb{E}L^m(t,0) \leq (m!)^p [\mathbb{E}L(t,0)]^m, \tag{4.4}$$

$$\mathbb{E}I_t^m \leq (m!)^{2p} [\mathbb{E}I_t]^m. \tag{4.5}$$

PROOF. This time we pick up the discrete case defined by (1.7). Similarly to (2.9) and (2.14),

$$L(t,0) = \frac{1}{(2\pi)^d} \int_{[-\pi,\pi]^d} d\lambda \prod_{j=1}^p \int_0^t e^{-i\lambda\cdot X_j(s)}\, ds, \tag{4.6}$$

$$\sum_{x\in\mathbb{Z}^d} L^2(t,x)\, dx = \frac{1}{(2\pi)^d} \int_{[-\pi,\pi]^d} \left| \prod_{j=1}^p \int_0^t e^{i\lambda\cdot X_j(s)}\, ds \right|^2 d\lambda. \tag{4.7}$$

Let $\Sigma_m$ be the group of permutations on $\{1,\dots,m\}$. By (4.6),

$$\begin{aligned}
&\mathbb{E}L^m(t,0)\\
&\quad = \frac{1}{(2\pi)^{md}} \int_{([-\pi,\pi]^d)^m} d\lambda_1\cdots d\lambda_m \left[\mathbb{E}\prod_{k=1}^m \int_0^t e^{i\lambda_k\cdot X(s)}\, ds\right]^p\\
&\quad = \frac{1}{(2\pi)^{md}} \int_{([-\pi,\pi]^d)^m} d\lambda_1\cdots d\lambda_m\\
&\qquad \times \left[\sum_{\sigma\in\Sigma_m} \int_{\{0\leq s_1\leq\cdots\leq s_m\leq t\}} \mathbb{E}\exp\left\{i\sum_{k=1}^m \lambda_{\sigma(k)}\cdot X(s_k)\right\} ds_1\cdots ds_m\right]^p\\
&\quad \leq \frac{(m!)^p}{(2\pi)^{md}} \int_{([-\pi,\pi]^d)^m} d\lambda_1\cdots d\lambda_m\\
&\qquad \times \left[\int_{\{0\leq s_1\leq\cdots\leq s_m\leq t\}} \mathbb{E}\exp\left\{i\sum_{k=1}^m \lambda_k\cdot X(s_k)\right\} ds_1\cdots ds_m\right]^p\\
&\quad = \frac{(m!)^p}{(2\pi)^{md}} \int_{([-\pi,\pi]^d)^m} d\lambda_1\cdots d\lambda_m\\
&\qquad \times \left[\int_{\{0\leq s_1\leq\cdots\leq s_m\leq t\}} \mathbb{E}\exp\left\{i\sum_{k=1}^m \left(\sum_{j=k}^m \lambda_j\right)\right.\right.\\
&\qquad\qquad \left.\left. \cdot (X(s_k)-X(s_{k-1}))\right\} ds_1\cdots ds_m\right]^p\\
&\quad = \frac{(m!)^p}{(2\pi)^{md}} \int_{([-\pi,\pi]^d)^m} d\lambda_1\cdots d\lambda_m
\end{aligned}$$

$$\times \left[\int_{\{0\le s_1\le\cdots\le s_m\le t\}} \exp\left\{-\sum_{k=1}^m (s_k - s_{k-1})\psi\left(\sum_{j=k}^m \lambda_j\right)\right\} ds_1\cdots ds_m\right]^p,$$

where the third step follows from the Jensen inequality and the index permutation and where we adopt the convention that $s_0 = 0$.

Write

$$f(\lambda_1,\dots,\lambda_m) = \left[\int_{\{0\le s_1\le\cdots\le s_m\le t\}} \exp\left\{-\sum_{k=1}^m (s_k - s_{k-1})\psi(\lambda_k)\right\} ds_1\cdots ds_m\right]^p,$$

where $\lambda_1,\dots,\lambda_m \in \mathbb{R}^d$. By (1.7), we have that, for any $\mathbf{k}_1,\dots,\mathbf{k}_m \in \mathbb{Z}^d$ and for any $\lambda_1,\dots,\lambda_m \in \mathbb{R}^d$,

$$f(\lambda_1 + (2\pi)\mathbf{k}_1,\dots,\lambda_m + (2\pi)\mathbf{k}_m) = f(\lambda_1,\dots,\lambda_m). \tag{4.8}$$

Notice that

$$\begin{aligned}\int_{([-\pi,\pi]^d)^m} d\lambda_1\cdots d\lambda_m f\left(\sum_{k=1}^m \lambda_k, \sum_{k=2}^m \lambda_k,\dots,\lambda_m\right)\\ = \int_{([-\pi,\pi]^d)^{m-1}} d\lambda_2\cdots d\lambda_m \int_{[-\pi,\pi]^d} f\left(\sum_{k=1}^m \lambda_k, \sum_{k=2}^m \lambda_k,\cdots,\lambda_m\right) d\lambda_1.\end{aligned}$$

By periodicity,

$$\int_{[-\pi,\pi]^d} f\left(\sum_{k=1}^m \lambda_k, \sum_{k=2}^m \lambda_k,\dots,\lambda_m\right) d\lambda_1 = \int_{[-\pi,\pi]^d} f\left(\lambda_1, \sum_{k=2}^m \lambda_k,\dots,\lambda_m\right) d\lambda_1.$$

So we have

$$\begin{aligned}\int_{([-\pi,\pi]^d)^m} d\lambda_1\cdots d\lambda_m f\left(\sum_{k=1}^m \lambda_k, \sum_{k=2}^m \lambda_k,\dots,\lambda_m\right)\\ = \int_{([-\pi,\pi]^d)^m} d\lambda_1\cdots d\lambda_m\, f\left(\lambda_1, \sum_{k=2}^m \lambda_k,\dots,\lambda_m\right).\end{aligned}$$

Repeating this procedure, we obtain

$$\begin{aligned}\int_{([-\pi,\pi]^d)^m} d\lambda_1\cdots d\lambda_m\, f\left(\sum_{k=1}^m \lambda_k, \sum_{k=2}^m \lambda_k,\dots,\lambda_m\right)\\ = \int_{([-\pi,\pi]^d)^m} d\lambda_1\cdots d\lambda_m\, f(\lambda_1,\dots,\lambda_m).\end{aligned} \tag{4.9}$$

In summary, we have proved that

$$\begin{aligned}
&\mathbb{E}L^m(t,0)\\
&\quad\le \frac{(m!)^p}{(2\pi)^{md}}\int_{([-\pi,\pi]^d)^m} d\lambda_1\cdots d\lambda_m\\
&\qquad\times\left[\int_{\{0\le s_1\le\cdots\le s_m\le t\}}\exp\left\{-\sum_{k=1}^m (s_k-s_{k-1})\psi(\lambda_k)\right\}ds_1\cdots ds_m\right]^p.
\end{aligned}$$

Notice that

$$\begin{aligned}
&\int_{\{0\le s_1\le\cdots\le s_m\le t\}}\exp\left\{-\sum_{k=1}^m (s_k-s_{k-1})\psi(\lambda_k)\right\}ds_1\cdots ds_m\\
&\quad=\int_{\substack{s_1+\cdots+s_m\le t\\ s_1,\dots,s_m\ge 0}}\exp\left\{-\sum_{k=1}^m s_k\psi(\lambda_k)\right\}ds_1\cdots ds_m\\
&\quad\le\int_0^t\cdots\int_0^t\exp\left\{-\sum_{k=1}^m s_k\psi(\lambda_k)\right\}ds_1\cdots ds_m\\
&\quad=\prod_{k=1}^m\int_0^t e^{-s\psi(\lambda_k)}\,ds.
\end{aligned}$$

Hence,

$$\begin{aligned}
\mathbb{E}L^m(t,0)&\le (m!)^p\left\{\frac{1}{(2\pi)^d}\int_{[-\pi,\pi]^d}d\lambda\left[\int_0^t e^{-s\psi(\lambda)}\,ds\right]^p\right\}^m\\
&=(m!)^p[\mathbb{E}L(t,0)]^m.
\end{aligned}$$

We now prove (4.5). By (4.7),

$$\begin{aligned}
&\mathbb{E}\left[\sum_{x\in\mathbb{Z}^d}L^2(t,x)\,dx\right]^m\\
&\quad=\frac{1}{(2\pi)^{md}}\int_{([-\pi,\pi]^d)^m}d\lambda_1\cdots d\lambda_m\\
&\qquad\times\left[\mathbb{E}\left|\sum_{\sigma\in\Sigma_m}\int_{\{0\le s_1\le\cdots\le s_m\le t\}}\exp\left\{i\sum_{k=1}^m\lambda_{\sigma(k)}\cdot X(s_k)\right\}ds_1\cdots ds_m\right|^2\right]^p\\
&\quad\le\frac{1}{(2\pi)^{md}}\int_{([-\pi,\pi]^d)^m}d\lambda_1\cdots d\lambda_m\\
&\qquad\times\Bigg[\sum_{\sigma\in\Sigma_m}\Bigg(\mathbb{E}\Bigg|\int_{\{0\le s_1\le\cdots\le s_m\le t\}}\exp\Bigg\{i\sum_{k=1}^m\lambda_{\sigma(k)}
\end{aligned}$$

$$\cdot X(s_k)\Bigg\} ds_1 \cdots ds_m\Bigg|^2\Bigg)^{1/2}\Bigg]^{2p}$$

$$\leq \frac{(m!)^{2p}}{(2\pi)^{md}} \int_{([-\pi,\pi]^d)^m} d\lambda_1 \cdots d\lambda_m$$

$$\times \left[\mathbb{E}\left|\int_{\{0\leq s_1\leq\cdots\leq s_m\leq t\}} \exp\left\{i\sum_{k=1}^m \lambda_k \cdot X(s_k)\right\} ds_1\cdots ds_m\right|^2\right]^p$$

$$= \frac{(m!)^{2p}}{(2\pi)^{md}} \int_{([-\pi,\pi]^d)^m} d\lambda_1 \cdots d\lambda_m$$

$$\times \Bigg[\mathbb{E}\Bigg|\int_{\{0\leq s_1\leq\cdots\leq s_m\leq t\}} \exp\Bigg\{i\sum_{k=1}^m \Bigg(\sum_{j=k}^m \lambda_j\Bigg)$$

$$\cdot (X(s_k) - X(s_{k-1}))\Bigg\} ds_1\cdots ds_m\Bigg|^2\Bigg]^p.$$

Since $X(t)$ takes only lattice values, the function

$$f(\lambda_1, \dots, \lambda_m)$$

$$\equiv \left[\mathbb{E}\left|\int_{\{0\leq s_1\leq\cdots\leq s_m\leq t\}} \exp\left\{i\sum_{k=1}^m \lambda_k \cdot (X(s_k) - X(s_{k-1}))\right\} ds_1\cdots ds_m\right|^2\right]^p$$

satisfies (4.8). Consequently, (4.9) holds. So we have

$$\mathbb{E}\left[\sum_{x\in\mathbb{Z}^d} L^2(t,x)\,dx\right]^m$$

$$\leq \frac{(m!)^{2p}}{(2\pi)^{md}} \int_{([-\pi,\pi]^d)^m} d\lambda_1 \cdots d\lambda_m$$

$$\times \Bigg[\mathbb{E}\Bigg|\int_{\{0\leq s_1\leq\cdots\leq s_m\leq t\}} \exp\Bigg\{i\sum_{k=1}^m \lambda_k \cdot (X(s_k)$$

$$- X(s_{k-1}))\Bigg\} ds_1\cdots ds_m\Bigg|^2\Bigg]^p.$$

Let $\widetilde{X}_1(t), \dots, \widetilde{X}_m(t)$ be independent copies of $X(t)$. Then for any $\lambda_1, \dots, \lambda_m \in \mathbb{R}^d$,

$$\int_{\{0\leq s_1\leq\cdots\leq s_m\leq t\}} \exp\left\{i\sum_{k=1}^m \lambda_k \cdot (X(s_k) - X(s_{k-1}))\right\} ds_1\cdots ds_m$$

$$\stackrel{d}{=} \int_{\{0\le s_1\le\cdots\le s_m\le t\}} \exp\left\{i\sum_{k=1}^m \lambda_k\cdot \widetilde{X}_k(s_k-s_{k-1})\right\} ds_1\cdots ds_m$$

$$= \int_{\substack{s_1+\cdots+s_m\le t\\ s_1,\dots,s_m\ge 0}} \exp\left\{i\sum_{k=1}^m \lambda_k\cdot \widetilde{X}_k(s_k)\right\} ds_1\cdots ds_m.$$

By Lemma 2.1, we have

$$\mathbb{E}\left|\int_{\{s_1+\cdots+s_m\le t\}} \exp\left\{i\sum_{k=1}^m \lambda_k\cdot \widetilde{X}_k(s_k)\right\} ds_1\cdots ds_m\right|^2$$

$$\le \mathbb{E}\left|\int_0^t\cdots\int_0^t \exp\left\{i\sum_{k=1}^m \lambda_k\cdot \widetilde{X}_k(s_k)\right\} ds_1\cdots ds_m\right|^2$$

$$= \prod_{k=1}^m \mathbb{E}\left|\int_0^t e^{i\lambda_k\cdot X(s)}\,ds\right|^2.$$

Summarizing what we have proved,

$$\mathbb{E}\left[\sum_{x\in\mathbb{Z}^d} L^2(t,x)\right]^m$$

$$\le (m!)^{2p}\left[\frac{1}{(2\pi)^d}\int_{[-\pi,\pi]^d}\mathbb{E}\left|\int_0^t e^{i\lambda_k\cdot X(s)}\,ds\right|^2 d\lambda\right]^m$$

$$= (m!)^{2p}\left[\mathbb{E}\sum_{x\in\mathbb{Z}^d} L^2(t,x)\right]^m.$$

□

LEMMA 4.2. *Under both (1.6) and (1.7), for any integers $a\ge 2$ and $m\ge 1$ and for any $t_1,\dots,t_a\ge 0$,*

$$[\mathbb{E}L^m(t_1+\cdots+t_a,0)]^{1/p} \le \sum_{\substack{k_1+\cdots+k_a=m\\ k_1,\dots,k_a\ge 0}} \frac{m!}{k_1!\cdots k_a!}[\mathbb{E}L^{k_1}(t_1,0)]^{1/p}\cdots[\mathbb{E}L^{k_a}(t_a,0)]^{1/p}, \tag{4.10}$$

$$[\mathbb{E}I^m_{t_1+\cdots+t_a}]^{1/(2p)} \le \sum_{\substack{k_1+\cdots+k_a=m\\ k_1,\dots,k_a\ge 0}} \frac{m!}{k_1!\cdots k_a!}[\mathbb{E}I^{k_1}_{t_1}]^{1/(2p)}\cdots[\mathbb{E}I^{k_a}_{t_a}]^{1/(2p)}, \tag{4.11}$$

*where $I_t$ is defined by (1.15).*

*Consequently, for any $\theta>0$,*

$$\sum_{m=0}^\infty \frac{\theta^m}{m!}[\mathbb{E}L^m(t_1+\cdots+t_a,0)]^{1/p} \le \prod_{k=1}^a\sum_{m=0}^\infty\frac{\theta^m}{m!}[\mathbb{E}L^m(t_k,0)]^{1/p}, \tag{4.12}$$

$$\sum_{m=0}^{\infty}\frac{\theta^m}{m!}[\mathbb{E}I^m_{t_1+\cdots+t_a}]^{1/(2p)}\le\prod_{k=1}^{a}\sum_{m=0}^{\infty}\frac{\theta^m}{m!}[\mathbb{E}I^m_{t_k}]^{1/(2p)}. \tag{4.13}$$

PROOF. Used inductively, only the establishment of (4.10) and (4.11) in the case $a=2$ is needed. That is, we only need to show that

$$[\mathbb{E}L^m(t_1+t_2,0)]^{1/p}\le\sum_{k=0}^{m}\binom{m}{k}[\mathbb{E}L^k(t_1,0)]^{1/p}[\mathbb{E}L^{m-k}(t_2,0)]^{1/p}, \tag{4.14}$$

$$[\mathbb{E}I^m_{t_1+t_2}]^{1/(2p)}\le\sum_{k=0}^{m}\binom{m}{k}[\mathbb{E}I^k_{t_1}]^{1/(2p)}[\mathbb{E}I^{m-k}_{t_2}]^{1/(2p)}. \tag{4.15}$$

We only consider the continuous case defined by (1.6), as the proof in the case (1.7) is analogous. Write

$$\Delta_1(\lambda)=\int_0^{t_1}e^{i\lambda\cdot X(t)}\,dt\quad\text{and}\quad\Delta_2(\lambda)=\int_{t_1}^{t_1+t_2}e^{i\lambda\cdot X(t)}\,dt,\qquad\lambda\in\mathbb{R}^d. \tag{4.16}$$

By (2.10),

$$\begin{aligned}&\mathbb{E}L^m(t_1+t_2,0)\\&\quad=\frac{1}{(2\pi)^{md}}\int_{(\mathbb{R}^d)^m}d\lambda_1\cdots d\lambda_m\left[\mathbb{E}\prod_{k=1}^{m}(\Delta_1(\lambda_k)+\Delta_2(\lambda_k))\right]^p\\&\quad=\frac{1}{(2\pi)^{md}}\int_{(\mathbb{R}^d)^m}d\lambda_1\cdots d\lambda_m\left[\sum_{l_1,\dots,l_m=1}^{2}\mathbb{E}(\Delta_{l_1}(\lambda_1)\cdots\Delta_{l_m}(\lambda_m))\right]^p.\end{aligned}$$

By Lemma 2.1, for any $(\lambda_1,\dots,\lambda_m)$ and $(l_1,\dots,l_m)$,

$$\mathbb{E}(\Delta_{l_1}(\lambda_1)\cdots\Delta_{l_m}(\lambda_m))\ge 0.$$

So the triangular inequality for the $L_p$-norm applies here, which gives

$$\begin{aligned}&[\mathbb{E}L^m(t_1+t_2,0)]^{1/p}\\&\quad\le\frac{1}{(2\pi)^{md}}\\&\qquad\times\sum_{l_1,\dots,l_m=1}^{2}\left\{\int_{(\mathbb{R}^d)^m}d\lambda_1\cdots d\lambda_m[\mathbb{E}(\Delta_{l_1}(\lambda_1)\cdots\Delta_{l_m}(\lambda_m))]^p\right\}^{1/p}.\end{aligned} \tag{4.17}$$

Let $(l_1,\dots,l_m)$ be arbitrary but fixed and let $k$ be the number of 1's among $l_1,\dots,l_m$:

$$\int_{(\mathbb{R}^d)^m}d\lambda_1\cdots d\lambda_m[\mathbb{E}(\Delta_{l_1}(\lambda_1)\cdots\Delta_{l_m}(\lambda_m))]^p$$

$$
\begin{aligned}
&= \int_{(\mathbb{R}^d)^m} d\lambda_1 \cdots d\lambda_m [\mathbb{E}(\Delta_1(\lambda_1)\cdots\Delta_1(\lambda_k))(\Delta_2(\lambda_{k+1})\cdots\Delta_2(\lambda_m))]^p \\
&= \int_{(\mathbb{R}^d)^m} d\lambda_1 \cdots d\lambda_m \Bigg[\mathbb{E}\Bigg(\Delta_1(\lambda_1)\cdots\Delta_1(\lambda_k) \\
&\qquad\qquad \times \exp\Bigg\{i\Bigg(\sum_{l=k+1}^{m}\lambda_l\Bigg)\cdot X(t_1)\Bigg\}\Bigg) \\
&\qquad\qquad \times \mathbb{E}(\widetilde{\Delta}(\lambda_{k+1})\cdots\widetilde{\Delta}(\lambda_m))\Bigg]^p \qquad (4.18)\\
&= \int_{(\mathbb{R}^d)^{m-k}} d\lambda_{k+1}\cdots d\lambda_m [\mathbb{E}(\widetilde{\Delta}(\lambda_{k+1})\cdots\widetilde{\Delta}(\lambda_m))]^p \\
&\quad \times \int_{(\mathbb{R}^d)^k} d\lambda_1\cdots d\lambda_k \Bigg[\mathbb{E}\Bigg(\Delta_1(\lambda_1)\cdots\Delta_1(\lambda_k) \\
&\qquad\qquad \times \exp\Bigg\{i\Bigg(\sum_{l=k+1}^{m}\lambda_l\Bigg)\cdot X(t_1)\Bigg\}\Bigg)\Bigg]^p,
\end{aligned}
$$

where

$$
\widetilde{\Delta}(\lambda) = \int_0^{t_2} e^{i\lambda\cdot X(s)}\,ds, \qquad \lambda \in \mathbb{R}^d, \tag{4.19}
$$

and where the second step follows from the independent increment property of Lévy processes.

For each $1 \le j \le p$, write

$$
\Delta_1^j(\lambda) = \int_0^{t_1} e^{i\lambda\cdot X_j(t)}\,ds, \qquad \lambda \in \mathbb{R}^d.
$$

Then for any $\lambda_{k+1},\dots,\lambda_m \in \mathbb{R}^d$,

$$
\begin{aligned}
&\int_{(\mathbb{R}^d)^k} d\lambda_1\cdots d\lambda_k \Bigg[\mathbb{E}\Bigg(\Delta_1(\lambda_1)\cdots\Delta_1(\lambda_k)\exp\Bigg\{i\Bigg(\sum_{l=k+1}^{m}\lambda_l\Bigg)\cdot X(t_1)\Bigg\}\Bigg)\Bigg]^p \\
&\quad = \int_{(\mathbb{R}^d)^k} d\lambda_1\cdots d\lambda_k\, \mathbb{E}\prod_{j=1}^{p}\Bigg(\Delta_1^j(\lambda_1)\cdots\Delta_1^j(\lambda_k)\exp\Bigg\{i\Bigg(\sum_{l=k+1}^{m}\lambda_l\Bigg)\cdot X_j(t_1)\Bigg\}\Bigg) \\
&\quad = \mathbb{E}\Bigg\{\exp\Bigg\{i\Bigg(\sum_{l=k+1}^{m}\lambda_l\Bigg)\cdot\Bigg(\sum_{j=1}^{p}X_j(t_1)\Bigg)\Bigg\} \\
&\qquad \times \int_{(\mathbb{R}^d)^k} d\lambda_1\cdots d\lambda_k \prod_{j=1}^{p}(\Delta_1^j(\lambda_1)\cdots\Delta_1^j(\lambda_k))\Bigg\}.
\end{aligned}
$$

By (2.9),

$$\int_{(\mathbb{R}^d)^k} d\lambda_1 \cdots d\lambda_k \prod_{j=1}^p (\Delta_1^j(\lambda_1)\cdots\Delta_1^j(\lambda_k))$$
$$= \left[\int_{\mathbb{R}^d} d\lambda \prod_{j=1}^p \Delta_1^j(\lambda)\right]^k = (2\pi)^{kd} L^k(t_1,0) \geq 0.$$

So we have

$$\int_{(\mathbb{R}^d)^k} d\lambda_1 \cdots d\lambda_k \left[\mathbb{E}\left(\Delta_1(\lambda_1)\cdots\Delta_1(\lambda_k)\exp\left\{i\left(\sum_{l=k+1}^m \lambda_l\right)\cdot X(t_1)\right\}\right)\right]^p$$
$$\leq (2\pi)^{kd}\mathbb{E}L^k(t_1,0).$$

We now come back to (4.18). By Lemma 2.1, for any $\lambda_{k+1},\dots,\lambda_m \in \mathbb{R}^d$,

$$[\mathbb{E}(\widetilde{\Delta}(\lambda_{k+1})\cdots\widetilde{\Delta}(\lambda_m))]^p \geq 0.$$

Therefore,

$$\begin{aligned}\int_{(\mathbb{R}^d)^m} & d\lambda_1 \cdots d\lambda_m [\mathbb{E}(\Delta_{l_1}(\lambda_1)\cdots\Delta_{l_m}(\lambda_m))]^p \\ &\leq (2\pi)^{kd}\mathbb{E}L^k(t_1,0)\int_{(\mathbb{R}^d)^{m-k}} d\lambda_{k+1}\cdots d\lambda_m \\ &\qquad \times [\mathbb{E}(\widetilde{\Delta}(\lambda_{k+1})\cdots\widetilde{\Delta}(\lambda_m))]^p \\ &= (2\pi)^{md}[\mathbb{E}L^k(t_1,0)][\mathbb{E}L^{m-k}(t_2,0)],\end{aligned} \tag{4.20}$$

where the last step follows from (2.10). Thus, (4.14) follows from (4.17) and (4.20).

We now come to the proof of (4.15). Recall that $\Delta_1(\lambda)$ and $\Delta_2(\lambda)$ are defined by (4.16). By (2.14),

$$\begin{aligned}\mathbb{E}&\left[\int_{\mathbb{R}^d} L^2(t_1+t_2,x)\,dx\right]^m \\ &= \frac{1}{(2\pi)^{md}}\int_{(\mathbb{R}^d)^m} d\lambda_1\cdots d\lambda_m \left[\mathbb{E}\left|\prod_{k=1}^m (\Delta_1(\lambda_k)+\Delta_2(\lambda_k))\right|^2\right]^p \\ &= \frac{1}{(2\pi)^{md}}\int_{(\mathbb{R}^d)^m} d\lambda_1\cdots d\lambda_m \left[\mathbb{E}\left|\sum_{l_1,\dots,l_m=1}^2 \Delta_{l_1}(\lambda_1)\cdots\Delta_{l_m}(\lambda_m)\right|^2\right]^p \\ &\leq \frac{1}{(2\pi)^{md}}\int_{(\mathbb{R}^d)^m} d\lambda_1\cdots d\lambda_m \\ &\qquad \times \left[\sum_{l_1,\dots,l_m=1}^2 \{\mathbb{E}(|\Delta_{l_1}(\lambda_1)|\cdots|\Delta_{l_m}(\lambda_m)|)^2\}^{1/2}\right]^{2p},\end{aligned}$$

where the last step follows from the triangular inequality for the $L_2$-norm. Using the triangular inequality for the $L_{2p}$-norm,

$$\begin{aligned}
&\left\{\mathbb{E}\left[\int_{\mathbb{R}^d} L^2(t_1+t_2,x)\,dx\right]^m\right\}^{1/2p} \\
&\quad\le \frac{1}{(2\pi)^{md}} \sum_{l_1,\dots,l_m=1}^{2} \left\{\int_{(\mathbb{R}^d)^m} d\lambda_1\cdots d\lambda_m \right. \\
&\qquad\qquad \left. \times [\mathbb{E}(|\Delta_{l_1}(\lambda_1)|\cdots|\Delta_{l_m}(\lambda_m)|)^2]^p\right\}^{1/2p}.
\end{aligned} \tag{4.21}$$

Let $(l_1,\dots,l_m)$ be arbitrary but fixed and let $k$ be the number of 1's among $l_1,\dots,l_m$:

$$\begin{aligned}
&\int_{(\mathbb{R}^d)^m} d\lambda_1\cdots d\lambda_m[\mathbb{E}(|\Delta_{l_1}(\lambda_1)|\cdots|\Delta_{l_m}(\lambda_m)|)^2]^p \\
&\quad= \int_{(\mathbb{R}^d)^m} d\lambda_1\cdots d\lambda_m[\mathbb{E}(|\Delta_1(\lambda_1)|\cdots|\Delta_1(\lambda_k)|)^2 \\
&\qquad\qquad \times (|\Delta_2(\lambda_{k+1})|\cdots|\Delta_2(\lambda_m)|)^2]^p.
\end{aligned}$$

By the identity

$$\begin{aligned}
|\Delta_2(\lambda)| &= \left|\int_0^{t_2} e^{i\lambda\cdot X(t_1)} e^{i\lambda\cdot(X(t_1+s)-X(t_1)}\,ds\right| \\
&= \left|\int_0^{t_2} e^{i\lambda\cdot(X(t_1+s)-X(t_1))}\,ds\right|,
\end{aligned}$$

we have that

$$\{|\Delta_2(\lambda)|;\lambda\in\mathbb{R}^d\} \overset{d}{=} \left\{\left|\int_0^{t_2} e^{i\lambda\cdot X(s)}\,ds\right|;\lambda\in\mathbb{R}^d\right\}$$

and that the two families

$$\{|\Delta_1(\lambda)|;\lambda\in\mathbb{R}^d\} \quad\text{and}\quad \{|\Delta_2(\lambda)|;\lambda\in\mathbb{R}^d\}$$

are independent. Hence,

$$\begin{aligned}
&\int_{(\mathbb{R}^d)^m} d\lambda_1\cdots d\lambda_m[\mathbb{E}(|\Delta_{l_1}(\lambda_1)|\cdots|\Delta_{l_m}(\lambda_m)|)^2]^p \\
&\quad= \left\{\int_{(\mathbb{R}^d)^k} d\lambda_1\cdots d\lambda_k\left[\mathbb{E}\left|\prod_{j=1}^k\int_0^{t_1} e^{i\lambda_j\cdot X(s)}\,dx\right|^2\right]^p\right\} \\
&\qquad \times \left\{\int_{(\mathbb{R}^d)^{m-k}} d\lambda_{k+1}\cdots d\lambda_m\left[\mathbb{E}\left|\prod_{j=k+1}^m\int_0^{t_2} e^{i\lambda_j\cdot X(s)}\,dx\right|^2\right]^p\right\}
\end{aligned} \tag{4.22}$$

$$= (2\pi)^{md}\mathbb{E}\left[\int_{\mathbb{R}^d} L^2(t_1,x)\,dx\right]^k \mathbb{E}\left[\int_{\mathbb{R}^d} L^2(t_2,x)\,dx\right]^{m-k}.$$

Finally, (4.5) follows from (4.21) and (4.22). □

The idea used in the above proof can be used to establish some similar results. Let $f(\lambda)\ge 0$ be bounded on $\mathbb{R}^d$ and write

$$L_t(f) = \begin{cases} \displaystyle\int_{\mathbb{R}^d} f(\lambda)\left[\prod_{j=1}^p \int_0^t e^{i\lambda\cdot X_j(s)}\,ds\right] d\lambda, & \text{under (1.6)},\\ \displaystyle\int_{[-\pi,\pi]^d} f(\lambda)\left[\prod_{j=1}^p \int_0^t e^{i\lambda\cdot X_j(s)}\,ds\right] d\lambda, & \text{under (1.7)}.\end{cases}$$

We state the following lemma without proof, as it is an obvious modification of the proof of (4.12).

LEMMA 4.3. *Under both (1.6) and (1.7), for any $t_1,\dots,t_a\ge 0$ and any integer $m\ge 1$,*

$$\begin{aligned}&[\mathbb{E}L^m_{t_1+\cdots+t_a}(f)]^{1/p}\\ &\quad\le \sum_{\substack{k_1+\cdots+k_a=m\\ k_1,\dots,k_a\ge 0}} \frac{m!}{k_1!\cdots k_a!}[\mathbb{E}|L_{t_1}(f)|^{k_1}]^{1/p}\cdots[\mathbb{E}|L_{t_a}(f)|^{k_a}]^{1/p}.\end{aligned} \tag{4.23}$$

*Consequently, for any $\theta>0$,*

$$\sum_{m=0}^\infty \frac{\theta^m}{m!}[\mathbb{E}L^m_{t_1+\cdots+t_a}(f)]^{1/p} \le \prod_{k=1}^a \sum_{m=0}^\infty \frac{\theta^m}{m!}[\mathbb{E}|L_{t_k}(f)|^m]^{1/p}. \tag{4.24}$$

REMARK. Inequalities (4.4) and (4.5) take a form similar to the inequality obtained in Lemma 3.1 of [38]; and (4.10) and (4.11) take a form similar to Theorem 5.1 of [7] and Theorem 6 of [8]. On the other hand, there are some differences at the technical level. First, all mentioned previous results are established in the setting of the random walks. Second, the proof of these results comes from a direct estimate of the local time (or intersection local times). In our setting, the estimate is carried out through Fourier transformation. Consequently, the property given in Lemma 2.1 is crucially needed. As a result, our argument can not be extended to the setting of additive random walks unless we put on some additional assumptions. For example, if we assume that

$$\mathbb{E}e^{i\lambda\cdot S(1)}\ge 0, \qquad \lambda\in\mathbb{R}^d, \tag{4.25}$$

then Lemmas 4.1 and 4.2 hold for the additive random walks. In the following we state a slightly different lemma in this direction.

To state the lemma, let $\{\omega_k\}_{k\ge 1}$ be an i.i.d. of symmetric random sequence with every finite moment and let $\{\omega_{1,k}\}_{k\ge 1},\dots,\{\omega_{p,k}\}_{k\ge 1}$ be its independent copies. We also assume independence between

$$\{\{S(k)\}_{k\ge 1};\{S_1(k)\}_{k\ge 1},\dots,\{S_p(k)\}_{k\ge 1}\}$$

and

$$\{\{\omega_k\}_{k\ge 1};\{\omega_{1,k}\}_{k\ge 1},\dots,\{\omega_{p,k}\}_{k\ge 1}\}.$$

Write

$$\xi(n,x)=\sum_{k_1,\dots,k_p=1}^{n}(\omega_{1,k_1}\cdots\omega_{p,k_p})\mathbb{1}_{\{S_1(k_1)+\cdots+S_p(k_p)=x\}},\qquad x\in\mathbb{Z}^d.$$

By symmetry, we have

$$\mathbb{E}\xi^{2m-1}(n,x)=0,\qquad m,n=1,2,\dots. \tag{4.26}$$

LEMMA 4.4. *Assume (4.25). For any integers $n_1,\dots,n_a\ge 1$ and $m\ge 1$,*

$$\begin{aligned}&[\mathbb{E}\xi^{2m}(n_1+\cdots+n_a,0)]^{1/p}\\ &\quad\le\sum_{\substack{k_1+\cdots+k_a=m\\ k_1,\dots,k_a\ge 0}}\frac{(2m)!}{(2k_1)!\cdots(2k_a)!}\\ &\qquad\qquad\times[\mathbb{E}\xi^{2k_1}(n_1,0)]^{1/p}\cdots[\mathbb{E}\xi^{2k_a}(n_a,0)]^{1/p}.\end{aligned}\tag{4.27}$$

*Consequently, for any $\theta>0$,*

$$\begin{aligned}&\sum_{m=0}^{\infty}\frac{\theta^m}{(2m)!}[\mathbb{E}\xi^{2m}(n_1+\cdots+n_a,0)]^{1/p}\\ &\qquad\le\prod_{k=1}^{a}\sum_{m=0}^{\infty}\frac{\theta^m}{(2m)!}[\mathbb{E}\xi^{2m}(n_k,0)]^{1/p}.\end{aligned}\tag{4.28}$$

PROOF. Notice that, for any $m\ge 1$ and $n\ge 1$,

$$\mathbb{E}\xi^{2m}(n,0)=\frac{1}{(2\pi)^{2md}}\int_{([-\pi,\pi]^d)^{2m}}d\lambda_1\cdots d\lambda_{2m}\left[\mathbb{E}\prod_{k=1}^{2m}\sum_{l=1}^{n}\omega_l e^{i\lambda\cdot S(l)}\right]^p. \tag{4.29}$$

To use the argument used in the proof of (4.14), we need the following property: for any $\lambda_1,\dots,\lambda_l\in\mathbb{R}^d$ and for any integer $j_1,\dots,j_l\ge 1$,

$$\mathbb{E}\prod_{k=1}^{l}\omega_{j_k}e^{i\lambda_k\cdot S(j_k)}\ge 0. \tag{4.30}$$

Indeed,

$$\mathbb{E}\prod_{k=1}^{l}\omega_{j_k}e^{i\lambda_k\cdot S(j_k)}=\mathbb{E}(\omega_{j_1}\cdots\omega_{j_l})\mathbb{E}\exp\left\{i\sum_{k=1}^{l}\lambda_k\cdot S(j_k)\right\}.$$

By (4.25) and by an argument similar to the one used in Lemma 2.1, we have

$$\mathbb{E}\exp\left\{i\sum_{k=1}^{l}\lambda_k\cdot S(j_k)\right\}\ge 0.$$

Write

$$\omega_{j_1}\cdots\omega_{j_l}=\omega_{k_1}^{l_1}\cdots\omega_{k_r}^{l_r},$$

where $k_1,\dots,k_r;l_1,\dots,l_r\ge 1$ are integers and $k_1,\dots,k_r$ are distinct. If any of $l_1,\dots,l_r$ are odd number, then by symmetry and independence, $\mathbb{E}(\omega_{j_1}\cdots\omega_{j_l})=0$. Otherwise we have $\mathbb{E}(\omega_{j_1}\cdots\omega_{j_l})\ge 0$.

Finally, by a proof similar to the proof of (4.14),

$$\begin{aligned}
&[\mathbb{E}\xi^{2m}(n_1+n_2,0)]^{1/p}\\
&\quad\le\sum_{k=0}^{2m}\binom{2m}{k}[\mathbb{E}|\xi(t_1,0)|^k]^{1/p}[\mathbb{E}\xi^{2m-k}(t_2,0)]^{1/p}\\
&\quad=\sum_{k=0}^{m}\binom{2m}{2k}[\mathbb{E}\xi^{2k}(t_1,0)]^{1/p}[\mathbb{E}\xi^{2(m-k)}(t_2,0)]^{1/p},
\end{aligned}$$

where the second step follows from (4.26). □

The following lemma is an application of Lemma 4.1 and Lemma 4.2.

Lemma 4.5. *In the assumptions of Theorem 1.1, for any $\delta>0$, there is a $C>0$ such that,*

$$\mathbb{E}L^m(t,0)\le C(m!)^p\delta^m t^{mp}a(t)^{-md},\tag{4.31}$$

$$\mathbb{E}I_t^m\le C(m!)^{2p}\delta^m t^{2mp}a(t)^{-md}\tag{4.32}$$

*for all $m\ge 1$ and sufficiently large $t$, where $I_t$ is defined by (1.15).*

*Consequently, for any $\theta>0$,*

$$\begin{aligned}
&\lim_{t\to\infty}\sum_{m=0}^{\infty}\frac{\theta^m}{m!}t^{-m}a(t)^{md/p}(\mathbb{E}L^m(t,0))^{1/p}\\
&\qquad=\sum_{m=0}^{\infty}\frac{\theta^m}{m!}(\mathbb{E}L_Y^m(1,0))^{1/p},
\end{aligned}\tag{4.33}$$

$$\lim_{t\to\infty}\sum_{m=0}^{\infty}\frac{\theta^m}{m!}t^{-m}a(t)^{md/2p}(\mathbb{E}I_t^m)^{1/2p}$$
$$= \sum_{m=0}^{\infty}\frac{\theta^m}{m!}\left\{\mathbb{E}\left[\int_{\mathbb{R}^d} L_Y^2(1,x)\,dx\right]^m\right\}^{1/(2p)}, \tag{4.34}$$

*where $L_Y(t,x)$ is the local time of the additive stable process generated by the stable process $Y$.*

PROOF. By the argument used in the proof of Theorem 3.1, one can see that there is $M>0$ such that

$$\mathbb{E}L(t,0)\le Mt^p a(t)^{-d} \quad \text{and} \quad \mathbb{E}I_t \le Mt^{2p}a(t)^{-d}$$

for large $t$. From (1.4) and (1.5), there is an integer $N\ge 1$ such that $N^{-p}a\times (N^{-1}t)^{-d}\le (2M)^{-1}\delta a(t)^{-d}$ as $t$ is large. By (4.10) in Lemma 4.2 and by (4.4) in Lemma 4.1,

$$\begin{aligned}
&[\mathbb{E}L^m(t,0)]^{1/p}\\
&\quad\le \sum_{\substack{k_1+\cdots+k_N=m\\ k_1,\dots,k_N\ge 0}}\frac{m!}{k_1!\cdots k_N!}[\mathbb{E}L^{k_1}(N^{-1}t,0)]^{1/p}\cdots[\mathbb{E}L^{k_N}(N^{-1}t,0)]^{1/p}\\
&\quad\le \sum_{\substack{k_1+\cdots+k_N=m\\ k_1,\dots,k_N\ge 0}}\frac{m!}{k_1!\cdots k_N!}k_1!(M(N^{-1}t)^p a(N^{-1}t)^{-d})^{k_1/p}\\
&\qquad\qquad\times\cdots\times k_N!(M(N^{-1}t)^p a(N^{-1}t)^{-d})^{k_N/p}\\
&\quad= m!(M(N^{-1}t)^p a(N^{-1}t)^{-d})^{m/p}\sum_{\substack{k_1+\cdots+k_N=m\\ k_1,\dots,k_N\ge 0}} 1\\
&\quad\le m!(2^{-1}\delta)^{m/p}t^m a(t)^{-md/p}\binom{m+N-1}{m}\\
&\quad\le C^{1/p}m!\delta^{m/p}t^m a(t)^{-md/p},
\end{aligned}$$

where the fourth step partially follows from the combinatorial fact that the equation $k_1+\cdots+k_N=m$ has $\binom{m+N-1}{m}$ $(\mathbb{Z}^+)^N$-valued solutions, and the last step follows from the generous bound

$$\binom{m+N-1}{m}\le C^{1/p}2^{m/p}, \qquad m=1,2,\dots.$$

The proof of (4.32) is similar. By the dominated convergence theorem and by Theorem 3.1, (4.33) and (4.34) follow from (4.31) and (4.32), respectively. □

We now prove (4.1). Let $s>0$ be fixed for a moment and let $n_t=[s^{-1}b_t]$, $r_t=t/n_t$. By (4.12) in Lemma 4.2,

$$\sum_{m=0}^\infty \frac{\theta^m}{m!}\left(\frac{b_t}{t}\right)^m a\left(\frac{t}{b_t}\right)^{md/p}(\mathbb{E}L^m(t,0))^{1/p} \\ \le \left\{\sum_{m=0}^\infty \frac{\theta^m}{m!}\left(\frac{b_t}{t}\right)^m a\left(\frac{t}{b_t}\right)^{md/p}(\mathbb{E}L^m(r_t,0))^{1/p}\right\}^{n_t}.$$

Take the regularity given in (1.4) into account:

$$\frac{b_t}{t}a\left(\frac{t}{b_t}\right)^{d/p}\sim s^{(\alpha p-d)/(\alpha p)}r_t^{-1}a(r_t)^{d/p}\qquad (t\to\infty).$$

By (4.33) in Lemma 4.5, therefore,

$$\sum_{m=0}^\infty \frac{\theta^m}{m!}\left(\frac{b_t}{t}\right)^m a\left(\frac{t}{b_t}\right)^{md/p}(\mathbb{E}L^m(r_t,0))^{1/p} \\ \longrightarrow \sum_{m=0}^\infty \frac{\theta^m}{m!}s^{(\alpha p-d)/(\alpha p)m}(\mathbb{E}L_Y^m(1,0))^{1/p}\qquad (t\to\infty).$$

Hence,

$$\begin{aligned}&\limsup_{t\to\infty}\frac{1}{b_t}\log\left\{\sum_{m=0}^\infty \frac{\theta^m}{m!}\left(\frac{b_t}{t}\right)^m a\left(\frac{t}{b_t}\right)^{md/p}(\mathbb{E}L^m(t,0))^{1/p}\right\}\\ &\qquad\le \frac{1}{s}\log\left\{\sum_{m=0}^\infty \frac{\theta^m}{m!}s^{(\alpha p-d)/(\alpha p)m}(\mathbb{E}L_Y^m(1,0))^{1/p}\right\}.\end{aligned}\tag{4.35}$$

In view of (1.10), by Theorem 5 in [8],

$$\begin{aligned}&\lim_{s\to\infty}\frac{1}{s}\log\left\{\sum_{m=0}^\infty \frac{\theta^m}{m!}s^{(\alpha p-d)/(\alpha p)m}(\mathbb{E}L_Y^m(1,0))^{1/p}\right\}\\ &\quad=\sup_{\lambda>0}\left\{\theta\lambda^{1/p}-p^{-1}\frac{d}{\alpha}(2\pi)^\alpha\left(1-\frac{d}{\alpha p}\right)^{(\alpha p-d)/d}\left(\frac{\lambda}{\rho_1}\right)^{\alpha/d}\right\}\\ &\quad=(2\pi)^{-\alpha d/(\alpha p-d)}\rho_1^{\alpha/(\alpha p-d)}\theta^{\alpha p/(\alpha p-d)}.\end{aligned}$$

Letting $s\to\infty$ on the right-hand side of (4.35) gives (4.1).

Similarly, applying (4.13) in Lemma 4.2, (4.34) in Lemma 4.5, Theorem 5 in [8] and (1.11), we can prove that

$$\begin{aligned}&\limsup_{t\to\infty}\frac{1}{b_t}\log\left\{\sum_{m=0}^\infty \frac{\theta^m}{m!}\left(\frac{b_t}{t}\right)^m a\left(\frac{t}{b_t}\right)^{md/2p}(\mathbb{E}I_t^m)^{1/2p}\right\}\\ &\qquad\le\sup_{\lambda>0}\left\{\theta\lambda^{1/2p}-(2p)^{-1}\frac{d}{2\alpha}(2\pi)^\alpha\left(1-\frac{d}{2\alpha p}\right)^{(2\alpha p-d)/d}\left(\frac{\lambda}{\rho_2}\right)^{\alpha/d}\right\}.\end{aligned}\tag{4.36}$$

We now claim that (4.36) implies (4.2). Indeed, by (ii) in Lemma 5.3 (with $p$ being replaced by $2p$) of [7], (4.36) leads to

$$\limsup_{t\to\infty}\frac{1}{b_t}\log\mathbb{P}\Big\{I_t\ge\lambda t^{2p}a\Big(\frac{t}{b_t}\Big)^{-d}\Big\}$$
$$\le-\frac{d}{2\alpha}(2\pi)^{\alpha}\Big(1-\frac{d}{2\alpha p}\Big)^{(2\alpha p-d)/d}\Big(\frac{\lambda}{\rho_2}\Big)^{\alpha/d} \tag{4.37}$$

or

$$\limsup_{t\to\infty}\frac{1}{b_t}\log\mathbb{P}\Big\{(I_t)^{1/2}\ge\lambda t^{p}a\Big(\frac{t}{b_t}\Big)^{-d/2}\Big\}$$
$$\le-(2\pi)^{\alpha}\frac{d}{2\alpha}\Big(1-\frac{d}{2\alpha p}\Big)^{(2\alpha p-d)/d}\Big(\frac{\lambda^2}{\rho_2}\Big)^{\alpha/d}. \tag{4.38}$$

Notice that

$$\sum_{m=0}^{\infty}\frac{\theta^m}{m!}\Big(\frac{b_t}{t}\Big)^m a\Big(\frac{t}{b_t}\Big)^{md/2p}\Big[\mathbb{E}I_t^{m/2}\mathbb{1}\Big\{(I_t)^{1/2}\ge\lambda t^p a\Big(\frac{t}{b_t}\Big)^{-d/2}\Big\}\Big]^{1/p}$$
$$\le\Big[\mathbb{P}\Big\{(I_t)^{1/2}\ge\lambda t^p a\Big(\frac{t}{b_t}\Big)^{-d/2}\Big\}\Big]^{1/2p}$$
$$\times\sum_{m=0}^{\infty}\frac{\theta^m}{m!}\Big(\frac{b_t}{t}\Big)^m a\Big(\frac{t}{b_t}\Big)^{md/2p}[\mathbb{E}I_t^m]^{1/2p}.$$

From (4.36) and (4.38),

$$\lim_{\lambda\to\infty}\lim_{t\to\infty}\frac{1}{b_t}\log\sum_{m=0}^{\infty}\frac{\theta^m}{m!}\Big(\frac{b_t}{t}\Big)^m a\Big(\frac{t}{b_t}\Big)^{md/2p}$$
$$\times\Big[\mathbb{E}I_t^{m/2}\mathbb{1}\Big\{(I_t)^{1/2}\ge\lambda t^p a\Big(\frac{t}{b_t}\Big)^{-d/2}\Big\}\Big]^{1/p}=-\infty.$$

By Lemma 5.3(i) of [8], therefore, (4.38) leads to

$$\limsup_{t\to\infty}\frac{1}{b_t}\log\sum_{m=0}^{\infty}\frac{\theta^m}{m!}\Big(\frac{b_t}{t}\Big)^m a\Big(\frac{t}{b_t}\Big)^{md/2p}\{\mathbb{E}I_t^{m/2}\}^{1/p}$$
$$\le\sup_{\lambda>0}\Big\{\theta\lambda^{1/p}-p^{-1}(2\pi)^{\alpha}\frac{d}{2\alpha}\Big(1-\frac{d}{2\alpha p}\Big)^{(2\alpha p-d)/d}\rho_2^{-\alpha/d}\lambda^{2\alpha/d}\Big\}$$
$$=(2\pi)^{-\alpha d/(2\alpha p-d)}\rho_2^{\alpha/(2\alpha p-d)}\theta^{2\alpha p/(2\alpha p-d)}.$$

We have proved (4.2).

Finally, we prove (4.3). Let $\{\tau_k\}_{k\ge0}$ be an i.i.d. sequence of exponential times such that $\mathbb{E}\tau_0=1$ and that $\{\tau_k\}_{k\ge0}$ is independent of $\{S(k)\}_{k\ge0}$. Write

$$T_0=0,\qquad T_k=\tau_0+\cdots+\tau_{k-1}$$

and define

$$N_t = k, \qquad \text{if } T_k \le t < T_{k+1}, \qquad k = 0, 1, 2, \dots.$$

It is well known that $N_t$ is a Poisson process with

$$\mathbb{P}\{N_t = n\} = e^{-t}\frac{t^n}{n!}, \qquad n = 0, 1, 2, \dots,$$

and the process $X(t)$ defined by

$$X(t) = S(N_t)$$

is a pure jump Lévy process with

$$\begin{aligned}\mathbb{E}\exp\{i\lambda \cdot X(t)\} &= \sum_{k=0} \mathbb{P}\{N_t = k\}\varphi(\lambda)^k \\ &= e^{-t}\sum_{k=0}\frac{t^k}{k!}\varphi(\lambda)^k = \exp\{-t[1-\varphi(\lambda)]\},\end{aligned}$$

where $\varphi(\lambda) = \mathbb{E}e^{i\lambda\cdot S(1)}$. In particular, the condition (1.7) is satisfied by $X(t)$. By (1.17) and by the classic law of large numbers for $N_t$, (1.4) holds.

Let $\{\tau_k^1\}_{k\ge 0}, \dots, \{\tau_k^p\}_{k\ge 0}$ be independent copies of $\{\tau_k\}_{k\ge 0}$. We assume the independence between

$$\{\{S(k)\}_{k\ge 0}; \{S_1(k)\}_{k\ge 0}, \dots, \{S_p(k)\}_{k\ge 0}\}$$

and

$$\{\{\tau_k\}_{k\ge 0}; \{\tau_k^1\}_{k\ge 0}, \dots, \{\tau_k^p\}_{k\ge 0}\}.$$

Let $(N_t^1, T_k^1), \dots, (N_t^p, T_k^p)$ be generated, respectively, by $\{\tau_k^1\}_{k\ge 0}, \dots, \{\tau_k^p\}_{k\ge 0}$ in the same way. Write $X_j(t) = S_j(N_t^j)$, $j = 1, \dots, p$.

From our construction,

$$\begin{aligned} L_n^x &\equiv \int_0^{T_{n+1}^1}\cdots\int_0^{T_{n+1}^p} \mathbb{1}_{\{X_1(s_1)+\cdots+X_p(s_p)=x\}}\,ds_1\cdots ds_p \\ &= \sum_{k_1,\dots,k_p=0}^{n}\int_{T_{k_1}^1}^{T_{k_1+1}^1}\cdots\int_{T_{k_p}^p}^{T_{k_p+1}^p}\mathbb{1}_{\{X_1(s_1)+\cdots+X_p(s_p)=x\}}\,ds_1\cdots ds_p \\ &= \sum_{k_1,\dots,k_p=0}^{n}(\tau_{k_1}^1\cdots\tau_{k_p}^p)\mathbb{1}_{\{S_1(k_1)+\cdots+S_p(k_p)=x\}}, \qquad x\in\mathbb{Z}^d. \end{aligned} \tag{4.39}$$

Let $\lambda > 1$ be fixed but arbitrary. By the Jensen inequality, for any $m, n \ge 1$,

$$\mathbb{E}\left\{\sum_{x\in\mathbb{Z}^d} L^2(\lambda n, x)\right\}^m$$

$$
\begin{aligned}
&= \mathbb{E}\Bigg\{\sum_{x\in\mathbb{Z}^d}\bigg[\int_0^{\lambda n}\cdots\int_0^{\lambda n}\mathbb{1}_{\{X_1(s_1)+\cdots+X_p(s_p)=x\}}\,ds_1\cdots ds_p\bigg]^2\Bigg\}^{m/2}\\
&\geq \mathbb{E}\Bigg\{\sum_{x\in\mathbb{Z}^d}\mathbb{1}_{\{\max_{1\leq j\leq p}T^j_{n+1}\leq\lambda n\}}\\
&\qquad\qquad \times\bigg[\int_0^{T^1_{n+1}}\cdots\int_0^{T^p_{n+1}}\mathbb{1}_{\{X_1(s_1)+\cdots+X_p(s_p)=x\}}\,ds_1\cdots ds_p\bigg]^2\Bigg\}^{m/2}\\
&\geq \mathbb{E}\Bigg\{\sum_{x\in\mathbb{Z}^d}\bigg[\sum_{k_1,\dots,k_p=0}^{n}\mathbb{E}(\tau^1_{k_1}\cdots\tau^p_{k_p})\mathbb{1}_{\{\max_{1\leq j\leq p}T^j_{n+1}\leq\lambda n\}}\\
&\qquad\qquad\qquad\qquad\qquad \times\mathbb{1}_{\{S_1(k_1)+\cdots+S_p(k_p)=x\}}\bigg]^2\Bigg\}^{m/2}\\
&= (\mathbb{E}\tau_1\mathbb{1}_{\{T_{n+1}\leq\lambda n\}})^{mp}\,\mathbb{E}\Bigg\{\sum_{x\in\mathbb{Z}^d}l^2(n,x)\Bigg\}^{m/2}.
\end{aligned}
$$

Notice that $\mathbb{E}\tau_1\mathbb{1}_{\{T_{n+1}\leq\lambda n\}}\to 1$ as $n\to\infty$. Applying (4.2) with $t=n\lambda$ and with $b_t$ being replaced by $\tilde{b}_t=b_{[\lambda^{-1}t]}$, and noticing

$$
\frac{b_n}{n}a\bigg(\frac{n}{b_n}\bigg)^{d/2p}\sim\lambda^{(2\alpha p-d)/(2\alpha p)}\frac{\tilde{b}_t}{t}a\bigg(\frac{t}{\tilde{b}_t}\bigg)^{d/2p}\qquad (n\to\infty),
$$

we have

$$
\begin{aligned}
&\limsup_{n\to\infty}\frac{1}{b_n}\log\sum_{m=0}^{\infty}\frac{\theta^m}{m!}\bigg(\frac{b_n}{n}\bigg)^m a\bigg(\frac{n}{b_n}\bigg)^{md/2p}\Bigg\{\mathbb{E}\bigg[\sum_{x\in\mathbb{Z}^d}l^2(n,x)\bigg]^{m/2}\Bigg\}^{1/p}\\
&\quad\leq\limsup_{n\to\infty}\frac{1}{b_n}\log\sum_{m=0}^{\infty}\frac{\theta^m}{m!}\bigg(\frac{b_n}{n}\bigg)^m a\bigg(\frac{n}{b_n}\bigg)^{md/2p}\Bigg\{\mathbb{E}\bigg[\sum_{x\in\mathbb{Z}^d}L^2(\lambda n,x)\bigg]^{m/2}\Bigg\}^{1/p}\\
&\quad\leq(2\pi)^{-\alpha d/(2\alpha p-d)}\rho_2^{\alpha/(2\alpha p-d)}\theta^{2\alpha p/(2\alpha p-d)}\lambda.
\end{aligned}
$$

Letting $\lambda\to 1^+$ on the right-hand side gives (4.3).

**5. Proof of (1.21), (1.22) and (1.24).** The proof of (1.22) and (1.24) is completely based on (1.14) and (1.20) and is the same as the proof of Theorem 1.2 in [9]. As for the upper bound of (1.21), it is a consequence of (1.13) of Theorem 1.1 in a standard practice of Borel–Cantelli lemma. By the argument used in the proof of Theorem 2.3 in [7], or Theorem 3 in [8],

to prove the lower bound

$$\limsup_{t\to\infty}\frac{1}{t^p}a\left(\frac{t}{\log\log t}\right)^d L(t,0) \tag{5.1}$$
$$\geq (2\pi)^{-d}\left(\frac{\alpha}{d}\right)^{d/\alpha}\left(1-\frac{d}{\alpha p}\right)^{-(p-d/\alpha)}\rho_1 \qquad \text{a.s.,}$$

we need only to establish the following lemma.

LEMMA 5.1. *Under the assumptions of Theorem 1.1, for any $\varepsilon>0$,*

$$\lim_{\delta\to 0^+}\limsup_{t\to\infty}\frac{1}{b_t}\log\sup_{|x|\leq\delta a(tb_t^{-1})}\mathbb{P}\left\{|L(t,0)-L(t,x)|\geq\varepsilon t^p a\left(\frac{t}{b_t}\right)^{-d}\right\} \tag{5.2}$$
$$=-\infty.$$

PROOF. We prove (5.2) under (1.7). Let

$$V(t,x)=\int_{[-\pi,\pi]^d}|1-e^{i\lambda\cdot x}|\left[\prod_{j=1}^p\int_0^t e^{i\lambda\cdot X_j(s)}\,ds\right]d\lambda, \qquad x\in\mathbb{Z}^d,\ t\geq 0.$$

Let $r_t=t/[b_t]$. By Lemma 4.3,

$$\sum_{m=0}^\infty\frac{\theta^m}{m!}\left(\frac{b_t}{t}\right)^m a\left(\frac{t}{b_t}\right)^{md/p}(\mathbb{E}V^m(t,x))^{1/p}$$
$$\leq\left\{\sum_{m=0}^\infty\frac{\theta^m}{m!}\left(\frac{b_t}{t}\right)^m a\left(\frac{t}{b_t}\right)^{md/p}(\mathbb{E}|V(r_t,x)|^m)^{1/p}\right\}^{[b_t]} \tag{5.3}$$
$$\leq\left\{\sum_{m=0}^\infty\frac{\theta^m}{m!}\left(\frac{b_t}{t}\right)^m a\left(\frac{t}{b_t}\right)^{md/p}(\mathbb{E}V^{2m}(r_t,x))^{1/2p}\right\}^{[b_t]}.$$

Notice that

$$\mathbb{E}V^{2m}(r_t,x)$$
$$=\int_{([-\pi,\pi]^d)^{2m}}d\lambda_1\cdots d\lambda_{2m}$$
$$\times\left(\prod_{k=1}^{2m}|1-e^{i\lambda_k\cdot x}|\right)\left[\mathbb{E}\prod_{k=1}^{2m}\int_0^{r_t}e^{i\lambda_k\cdot X(s)}\,ds\right]^p.$$

As $m=1$,

$$\mathbb{E}V^2(r_t,x)=a(r_t)^{2d}\int_{([-\pi,\pi]^d)^2}d\lambda_1\,d\lambda_2\left(\prod_{k=1}^2|1-e^{ia(r_t)^{-1}\lambda_k\cdot x}|\right)$$

$$\times \left[\mathbb{E}\prod_{k=1}^{2}\int_0^{r_t} e^{ia(r_t)^{-1}\lambda_k\cdot X(s)}\,ds\right]^p.$$

An argument similar to the one used for (3.5) gives

$$\limsup_{\delta\to 0^+}\limsup_{t\to\infty} t^{-2p} a\left(\frac{t}{b_t}\right)^{2d} \sup_{|x|\le \delta a(r_t)} \mathbb{E}V^2(r_t,x) = 0.$$

As for general $m\ge 1$,

$$\begin{aligned}\mathbb{E}V^{2m}(r_t,x) &\le 2^{2m}\mathbb{E}L^{2m}(r_t,0)\\ &\le C(m!)^{2p}\left(\frac{1}{2\theta}\right)^{2mp}\left(\frac{t}{b_t}\right)^{2mp} a\left(\frac{t}{b_t}\right)^{-2md}, \qquad x\in\mathbb{Z}^d,\end{aligned}$$

where the second step follows from (4.31) in Lemma 4.5. By the dominated convergence theorem,

$$\limsup_{\delta\to 0^+}\limsup_{t\to\infty}\sup_{|x|\le\delta a(r_t)}\sum_{m=0}^{\infty}\frac{\theta^m}{m!}\left(\frac{b_t}{t}\right)^m a\left(\frac{t}{b_t}\right)^{md/p}(\mathbb{E}V^{2m}(r_t,x))^{1/2p} = 1.$$

By (5.3), for any $\theta>0$,

$$\limsup_{\delta\to 0^+}\limsup_{t\to\infty}\frac{1}{b_t}\log\sup_{|x|\le\delta a(tb_t^{-1})}\sum_{m=0}^{\infty}\frac{\theta^m}{m!}\left(\frac{b_t}{t}\right)^m a\left(\frac{t}{b_t}\right)^{md/p}(\mathbb{E}V^m(t,x))^{1/p} = 0.$$

Consequently [notice that $\mathbb{E}V^m(t,x)\ge 0$ for all $m\ge 1$],

$$\begin{aligned}\limsup_{\delta\to 0^+}\limsup_{t\to\infty}\frac{1}{b_t}\log\sup_{|x|\le\delta a(tb_t^{-1})}\sum_{m=0}^{\infty}\frac{\theta^{2m}}{(2m)!}\left(\frac{b_t}{t}\right)^{2m} a\left(\frac{t}{b_t}\right)^{2md/p}&\\ \times(\mathbb{E}V^{2m}(t,x))^{1/p} = 0&\end{aligned}$$

On the other hand,

$$\begin{aligned}\mathbb{E}[L(t,0)-L(t,x)]^{2m} &= \frac{1}{(2\pi)^{2md}}\int_{([-\pi,\pi]^d)^{2m}} d\lambda_1\cdots d\lambda_{2m}\\ &\qquad\times\left(\prod_{k=1}^{2m}(1-e^{i\lambda_k\cdot x})\right)\left[\mathbb{E}\prod_{k=1}^{2m}\int_0^t e^{i\lambda_k\cdot X(s)}\,ds\right]^p\\ &\le \frac{1}{(2\pi)^{2md}}\mathbb{E}V^{2m}(t,x).\end{aligned}$$

Thus,

$$\begin{aligned}\lim_{\delta\to 0^+}\limsup_{t\to\infty}\frac{1}{b_t}\log\sup_{|x|\le\delta a(tb_t^{-1})}\sum_{m=0}^{\infty}\frac{\theta^{2m}}{(2m)!}\left(\frac{b_t}{t}\right)^{2m} a\left(\frac{t}{b_t}\right)^{2md/p}&\\ \times(\mathbb{E}[L(t,0)-L(t,x)]^{2m})^{1/p} = 0.&\end{aligned}$$

By the Chebyshev inequality,

$$(\varepsilon^{1/p}b_t)^{2m}\left(\mathbb{P}\left\{|L(t,0)-L(t,x)|\ge \varepsilon t^p a\left(\frac{t}{b_t}\right)^{-d}\right\}\right)^{1/p} \le \left(\frac{b_t}{t}\right)^{2m} a\left(\frac{t}{b_t}\right)^{2md/p}(\mathbb{E}[L(t,0)-L(t,x)]^{2m})^{1/p}. \tag{5.4}$$

Hence,

$$\left\{\sum_{m=0}^\infty \frac{(\lambda b_t)^{2m}}{(2m)!}\right\}\left(\mathbb{P}\left\{|L(t,0)-L(t,x)|\ge \varepsilon t^p a\left(\frac{t}{b_t}\right)^{-d}\right\}\right)^{1/p} \le \sum_{m=0}^\infty \frac{(\lambda\varepsilon^{-/p})^{2m}}{(2m)!}\left(\frac{b_t}{t}\right)^{2m} a\left(\frac{t}{b_t}\right)^{2md/p}(\mathbb{E}[L(t,0)-L(t,x)]^{2m})^{1/p}. \tag{5.5}$$

Therefore,

$$\lim_{\delta\to 0^+}\limsup_{t\to\infty}\frac{1}{b_t}\log\sup_{|x|\le\delta a(tb_t^{-1})}\mathbb{P}\left\{|L(t,0)-L(t,x)|\ge \varepsilon t^p a\left(\frac{t}{b_t}\right)^{-d}\right\}\le -p\lambda.$$

Letting $\lambda\to\infty$ on the right-hand side gives (5.2). □

**6. Proof of (1.19) and (1.23).** Let $\tau_k, N_t, T_k$ be defined as in Section 4 and recall that $X(t)=S(N_t)$ is a $\mathbb{Z}^d$-valued Lévy process satisfying the conditions given in Theorem 1.1. By Cramér's large deviation principle (Theorem 2.2.30, [18]), for any $\delta>0$, there is a $u>0$ such that

$$\mathbb{P}\{|T_{n+1}-n|\ge\delta n\}\le e^{-un}, \tag{6.1}$$

as $n$ is large. By the classical law of large numbers, $T_{n+1}\sim n$ a.s. as $n\to\infty$. Therefore, (1.13) in Theorem 1.1 and (1.21) in Theorem 1.3 imply, respectively,

$$\lim_{n\to\infty}\frac{1}{b_n}\log\mathbb{P}\left\{L_n^0\ge\lambda n^p a\left(\frac{n}{b_n}\right)^{-d}\right\} = -(2\pi)^\alpha\frac{d}{\alpha}\left(1-\frac{d}{\alpha p}\right)^{(\alpha p-d)/d}\left(\frac{\lambda}{\rho_1}\right)^{\alpha/d}, \tag{6.2}$$

$$\limsup_{n\to\infty}\frac{1}{n^p}a\left(\frac{n}{\log\log n}\right)^d L_n^0 = (2\pi)^{-d}\left(\frac{\alpha}{d}\right)^{d/\alpha}\left(1-\frac{d}{\alpha p}\right)^{-(p-d/\alpha)}\rho_1 \qquad \text{a.s.}, \tag{6.3}$$

where $L_n^x$ is defined in (4.39).

The only thing we need to do in this section is to show that, for any $\varepsilon > 0$,

$$\lim_{n\to\infty} \frac{1}{b_n} \log \mathbb{P}\left\{ |l(n,0) - L_n^0| \geq \varepsilon \lambda n^p a\left(\frac{n}{b_n}\right)^{-d} \right\} = -\infty, \tag{6.4}$$

which passes the moderate deviation and the law of the iterated logarithm from $L_n^0$ to $l(n,0)$, and therefore leads to (1.19) and (1.23).

By (4.39),

$$\begin{aligned} |L_n^0 - l(n,0)| &\leq \sum_{j=1}^p \left| \sum_{k_1,\dots,k_p=0}^n (\tau_{k_j}^j - 1)(\tau_{k_{j+1}}^{j+1} \cdots \tau_{k_p}^p) \mathbb{1}_{\{S_1(k_1)+\cdots+S_p(k_p)=0\}} \right| \\ &= \sum_{j=1}^p J_j(n), \qquad \text{say.} \end{aligned}$$

By the triangular inequality and by an estimate similar to the one carried out in (5.4) and (5.5), we will have (6.4) if we prove that, for any $1 \leq j \leq p$ and any $\theta > 0$,

$$\lim_{n\to\infty} \frac{1}{b_n} \log \sum_{m=0}^\infty \frac{\theta^{2m}}{(2m)!} \left(\frac{b_n}{n}\right)^{2m} a\left(\frac{n}{b_n}\right)^{2md/p} (\mathbb{E} J_j^{2m}(n))^{1/p} = 0. \tag{6.5}$$

Let $\{\tau_k'\}_{k\geq 0}$ be a copy of $\{\tau_k\}_{k\geq 0}$ and let $\{\varepsilon_k\}_{k\geq 0}$ be an i.i.d. sequence such that $\{\varepsilon_k = -1\} = \mathbb{P}\{\varepsilon_k = 1\} = 1/2$. We assume independence among all sequences we defined so far. By the Jensen inequality,

$$\begin{aligned} &\mathbb{E} J_j^{2m}(n) \\ &\quad \leq \mathbb{E}\left[ \sum_{k_1,\dots,k_p=0}^n (\tau_{k_j}^j - \tau_{k_j}')(\tau_{k_{j+1}}^{j+1} \cdots \tau_{k_p}^p) \mathbb{1}_{\{S_1(k_1)+\cdots+S_p(k_p)=0\}} \right]^{2m} \\ &\quad = \mathbb{E}\left[ \sum_{k_1,\dots,k_p=0}^n \varepsilon_{k_j} (\tau_{k_j}^j - \tau_{k_j}')(\tau_{k_{j+1}}^{j+1} \cdots \tau_{k_p}^p) \mathbb{1}_{\{S_1(k_1)+\cdots+S_p(k_p)=0\}} \right]^{2m} \\ &\quad \leq \left\{ \left( \mathbb{E}\left[ \sum_{k_1,\dots,k_p=0}^n (\varepsilon_{k_j} \tau_{k_j}^j)(\tau_{k_{j+1}}^{j+1} \cdots \tau_{k_p}^p) \mathbb{1}_{\{S_1(k_1)+\cdots+S_p(s_k)=0\}} \right]^{2m} \right)^{1/(2m)} \right. \\ &\qquad \left. + \left( \mathbb{E}\left[ \sum_{k_1,\dots,k_p=0}^n (\varepsilon_{k_j} \tau_{k_j}')(\tau_{k_{j+1}}^{j+1} \cdots \tau_{k_p}^p) \mathbb{1}_{\{S_1(k_1)+\cdots+S_p(k_p)=0\}} \right]^{2m} \right)^{1/(2m)} \right\}^{2m} \\ &\quad = 2^{2m} \mathbb{E}\left[ \sum_{k_1,\dots,k_p=0}^n (\varepsilon_{k_j} \tau_{k_j}^j)(\tau_{k_{j+1}}^{j+1} \cdots \tau_{k_p}^p) \mathbb{1}_{\{S_1(k_1)+\cdots+S_p(k_p)=0\}} \right]^{2m} \end{aligned}$$

$$\le 2^{2m}\mathbb{E}\left[\sum_{k_1,\ldots,k_p=0}^{n}(\tau_{k_1}^1\cdots\tau_{k_{j-1}}^{j-1})(\varepsilon_{k_j}\tau_{k_j}^{j})(\tau_{k_{j+1}}^{j+1}\cdots\tau_{k_p}^{p})\mathbb{1}_{\{S_1(k_1)+\cdots+S_p(k_p)=0\}}\right]^{2m}$$

$$= 2^{2m}\mathbb{E}\left[\sum_{k_1,\ldots,k_p=0}^{n}(\varepsilon_{k_j}\tau_{k_1}^1)(\tau_{k_2}^2\cdots\tau_{k_p}^{p})\mathbb{1}_{\{S_1(k_1)+\cdots+S_p(k_p)=0\}}\right]^{2m}.$$

Similarly to (4.39),

$$\sum_{k_1,\ldots,k_p=0}^{n}(\varepsilon_{k_j}\tau_{k_1}^1)(\tau_{k_2}^2\cdots\tau_{k_p}^{p})\mathbb{1}_{\{S_1(k_1)+\cdots+S_p(k_p)=0\}}$$
$$= \sum_{k=0}^{n}\varepsilon_k\tau_k^1\int_0^{T_{n+1}^2}\cdots\int_0^{T_{n+1}^p}\mathbb{1}_{\{S_1(k)+X_2(s_2)+\cdots+X_p(s_p)=0\}}\,ds_2\cdots ds_p.$$

Let $\gamma>1$ and use the notation "$\mathbb{E}^{\varepsilon}$" for the expectation with respect to the sequence $\{\varepsilon_k\}_{k\ge 0}$. By the contraction principle (see, e.g., Theorem 4.4, page 95, in [37]), as $T_{n+1}^2,\ldots,T_{n+1}^p\le\gamma n$,

$$\mathbb{E}^{\varepsilon}\left[\sum_{k=0}^{n}\varepsilon_k\tau_k^1\int_0^{T_{n+1}^2}\cdots\int_0^{T_{n+1}^p}\mathbb{1}_{\{S_1(k)+X_2(s_2)+\cdots+X_p(s_p)=0\}}\,ds_2\cdots ds_p\right]^{2m}$$
$$\le \mathbb{E}^{\varepsilon}\left[\sum_{k=0}^{n}\varepsilon_k\tau_k^1\int_0^{\gamma n}\cdots\int_0^{\gamma n}\mathbb{1}_{\{S_1(k)+X_2(s_2)+\cdots+X_p(s_p)=0\}}\,ds_2\cdots ds_p\right]^{2m}.$$

Hence,

$$\mathbb{E}\left[\sum_{k=1}^{n}\varepsilon_k\tau_k^1\int_0^{T_{n+1}^2}\cdots\int_0^{T_{n+1}^p}\mathbb{1}_{\{S_1(k)+X_2(s_2)+\cdots+X_p(s_p)=0\}}\,ds_2\cdots ds_p\right]^{2m}$$
$$\le 2^{2m}\mathbb{E}\Bigg[\sum_{k=0}^{n}\varepsilon_k\tau_k^1$$
$$\times\int_0^{\gamma n}\cdots\int_0^{\gamma n}\mathbb{1}_{\{S_1(k)+X_2(s_2)+\cdots+X_p(s_p)=0\}}\,ds_2\cdots ds_p\Bigg]^{2m}$$
$$+\mathbb{E}\Bigg\{\Bigg[\sum_{k=0}^{n}\varepsilon_k\tau_k^1$$
$$\times\int_0^{T_{n+1}^2}\cdots\int_0^{T_{n+1}^p}\mathbb{1}_{\{S_1(k)+X_2(s_2)+\cdots+X_p(s_p)=0\}}\,ds_2\cdots ds_p\Bigg]^{2m}$$
$$\times\mathbb{1}_{\{\max_{2\le j\le p}T_{n+1}^j\ge\gamma n\}}\Bigg\}.$$

Notice that the second term on the right-hand side is bounded by

$$\begin{aligned}&\mathbb{E}[(T_{n+1}^1\cdots T_{n+1}^p)^{2m}\mathbb{1}_{\{\max_{2\le j\le p}T_{n+1}^j\ge\gamma n\}}]\\&\qquad\le\mathbb{E}(T_{n+1})^{2m}[\mathbb{E}(T_{n+1})^{4m}]^{(p-1)/2}(p\mathbb{P}\{T_{n+1}\ge\gamma n\})^{1/2}.\end{aligned}$$

Write

$$K_n=\sum_{k=0}^n\varepsilon_k\tau_k^1\int_0^n\cdots\int_0^n\mathbb{1}_{\{S_1(k)+X_2(s_2)+\cdots+X_p(s_p)=0\}}\,ds_2\cdots ds_p.$$

In view of (6.1), it suffices to show that, for any $\theta>0$,

$$\lim_{n\to\infty}\frac{1}{b_n}\log\sum_{m=0}^\infty\frac{\theta^{2m}}{(2m)!}\left(\frac{b_n}{n}\right)^{2m}a\left(\frac{n}{b_n}\right)^{2md/p}(\mathbb{E}K_n^{2m})^{1/p}=0. \tag{6.6}$$

Using the Lévy inequality (Proposition 2.3, page 47, [37]) conditionally,

$$\begin{aligned}&\mathbb{E}K_n^{2m}\\&\quad\le 2^{2m}\mathbb{E}\left[\sum_{k=0}^{2n+1}\varepsilon_k\tau_k^1\int_0^n\cdots\int_0^n\mathbb{1}_{\{S_1(k)+X_2(s_2)+\cdots+X_p(s_p)=0\}}\,ds_2\cdots ds_p\right]^{2m}\\&\quad\le 4^{2m}\left\{\mathbb{E}\left[\sum_{k=0}^n\varepsilon_{2k}\tau_{2k}^1\int_0^n\cdots\int_0^n\mathbb{1}_{\{S_1(2k)+X_2(s_2)+\cdots+X_p(s_p)=0\}}\,ds_2\cdots ds_p\right]^{2m}\right.\\&\qquad\quad+\mathbb{E}\left[\sum_{k=0}^n\varepsilon_{2k+1}\tau_{2k+1}^1\right.\\&\qquad\qquad\left.\left.\times\int_0^n\cdots\int_0^n\mathbb{1}_{\{S_1(2k+1)+X_2(s_2)+\cdots+X_p(s_p)=0\}}\,ds_2\cdots ds_p\right]^{2m}\right\}.\end{aligned}$$

Notice that

$$\begin{aligned}&\mathbb{E}\left[\sum_{k=0}^n\varepsilon_{2k}\tau_{2k}^1\int_0^n\cdots\int_0^n\mathbb{1}_{\{S_1(2k)+X_2(s_2)+\cdots+X_p(s_p)=0\}}\,ds_2\cdots ds_p\right]^{2m}\\&\quad=(2\pi)^{-2md}\int_{([-\pi,\pi]^d)^{2m}}d\lambda_1\cdots d\lambda_{2m}\\&\qquad\times\left[\mathbb{E}\prod_{k=1}^{2m}\sum_{l=0}^n\varepsilon_l e^{i\lambda_k\cdot S(2l)}\right]\left[\mathbb{E}\prod_{k=1}^{2m}\int_0^n e^{i\lambda_k\cdot X(s)}\,ds\right]^{p-1}.\end{aligned}$$

Write $\tilde{S}(l)=S(2l)$. Then $\{\tilde{S}(l)\}_{l\ge0}$ is a random walk satisfying the condition (4.25), and therefore (4.30) with $\omega_l=\varepsilon_l\tau_l$. By Lemma 2.1 and the

Hölder inequality,

$$\mathbb{E}\left[\sum_{k=0}^{n}\varepsilon_{2k}\tau_{2k}^{1}\int_0^n\cdots\int_0^n \mathbb{1}_{\{S_1(2k)+X_2(s_2)+\cdots+X_p(s_p)=0\}}\,ds_2\cdots ds_p\right]^{2m}$$

$$\le (2\pi)^{-2md}\left\{\int_{([-\pi,\pi]^d)^{2m}} d\lambda_1\cdots d\lambda_{2m}\left[\mathbb{E}\prod_{k=1}^{2m}\sum_{l=0}^{n}\varepsilon_l e^{i\lambda_k\cdot\tilde{S}(l)}\right]^p\right\}^{1/p}$$

$$\times\left\{\int_{([-\pi,\pi]^d)^{2m}} d\lambda_1\cdots d\lambda_{2m}\left[\mathbb{E}\prod_{k=1}^{2m}\int_0^n e^{i\lambda_k\cdot X(s)}\,ds\right]^p\right\}^{(p-1)/p}$$

$$= (\mathbb{E}\tilde{\xi}^{2m}(n,0))^{1/p}(\mathbb{E}L^{2m}(n,0))^{(p-1)/p},$$

where

$$\tilde{\xi}(n,x)=\sum_{k_1,\dots,k_p=0}^{n}(\varepsilon_{k_1}^1\tau_k^1)\cdots(\varepsilon_{k_p}^p\tau_k^p)\mathbb{1}_{\{\tilde{S}_1(k_1)+\cdots+\tilde{S}_p(k_p)=x\}},\qquad x\in\mathbb{Z}^d,$$

and where $\{\varepsilon_k^1\}_{k\ge 0},\dots,\{\varepsilon_k^p\}_{k\ge 0}$ are independent copies of $\{\varepsilon_k\}_{k\ge 0}$.

Similarly,

$$\mathbb{E}\left[\sum_{k=0}^{n}\varepsilon_{2k+1}\tau_{2k-1}^{1}\int_0^n\cdots\int_0^n \mathbb{1}_{\{S_1(2k+1)+X_2(s_2)+\cdots+X_p(s_p)=0\}}\,ds_2\cdots ds_p\right]^{2m}$$

$$=(2\pi)^{-2md}\int_{([-\pi,\pi]^d)^{2m}} d\lambda_1\cdots d\lambda_{2m}$$

$$\times\left[\mathbb{E}\prod_{k=1}^{2m}\sum_{l=0}^{n}\varepsilon_l e^{i\lambda_k\cdot S(2l+1)}\right]\left[\mathbb{E}\prod_{k=1}^{2m}\int_0^n e^{i\lambda_k\cdot X(s)}\,ds\right]^{p-1}.$$

In view of (4.30) [with $S(k)$ being replaced by $\tilde{S}(k)$],

$$\left|\mathbb{E}\prod_{k=1}^{2m}\sum_{l=0}^{n}\varepsilon_l e^{i\lambda_k\cdot S(2l+1)}\right| = \left|\mathbb{E}e^{i(\lambda_1+\cdots+\lambda_m)\cdot S(1)}\mathbb{E}\prod_{k=1}^{2m}\sum_{l=0}^{n}\varepsilon_l e^{i\lambda_k\cdot\tilde{S}(l)}\right|$$

$$\le \mathbb{E}\prod_{k=1}^{2m}\sum_{l=0}^{n}\varepsilon_l e^{i\lambda_k\cdot\tilde{S}(l)}.$$

Therefore, we also have

$$\mathbb{E}\left[\sum_{k=0}^{n}\varepsilon_{2k+1}\tau_{2k+1}^{1}\int_0^n\cdots\int_0^n \mathbb{1}_{\{S_1(2k+1)+X_2(s_2)+\cdots+X_p(s_p)=0\}}\,ds_2\cdots ds_p\right]^{2m}$$

$$\le (\mathbb{E}\tilde{\xi}^{2m}(n,0))^{1/p}(\mathbb{E}L^{2m}(n,0))^{(p-1)/p}.$$

Summarizing the argument since (6.6), it remains to show that, for any $\theta > 0$,

$$\lim_{n\to\infty}\frac{1}{b_n}\log\sum_{m=0}^{\infty}\frac{\theta^{2m}}{(2m)!}\left(\frac{b_n}{n}\right)^{2m}a\left(\frac{n}{b_n}\right)^{2md/p} \times (\mathbb{E}\tilde{\xi}^{2m}(n,0))^{1/p^2}(\mathbb{E}L^{2m}(n,0))^{(p-1)/p^2}=0. \tag{6.7}$$

Let $q > 1$ be the conjugate number of $p$ defined by the relation $p^{-1}+q^{-1}=1$ and let $\delta > 0$ be fixed but arbitrary. By the Hölder inequality,

$$\begin{aligned}&\sum_{m=0}^{\infty}\frac{\theta^{2m}}{(2m)!}\left(\frac{b_n}{n}\right)^{2m}a\left(\frac{n}{b_n}\right)^{2md/p}(\mathbb{E}\tilde{\xi}^{2m}(n,0))^{1/p^2}(\mathbb{E}L^{2m}(n,0))^{(p-1)/p^2}\\&\quad\le\left\{\sum_{m=0}^{\infty}\frac{(\delta^{-1}\theta)^{2pm}}{(2m)!}\left(\frac{b_n}{n}\right)^{2m}a\left(\frac{n}{b_n}\right)^{2md/p}(\mathbb{E}\tilde{\xi}^{2m}(n,0))^{1/p}\right\}^{1/p}\\&\qquad\times\left\{\sum_{m=0}^{\infty}\frac{\delta^{2qm}}{(2m)!}\left(\frac{b_n}{n}\right)^{2m}a\left(\frac{n}{b_n}\right)^{2md/p}(\mathbb{E}L^{2m}(n,0))^{1/p}\right\}^{(p-1)/p}.\end{aligned}$$

Since

$$\begin{aligned}&\sum_{m=0}^{\infty}\frac{\delta^{2qm}}{(2m)!}\left(\frac{b_n}{n}\right)^{2m}a\left(\frac{n}{b_n}\right)^{2md/p}(\mathbb{E}L^{2m}(n,0))^{1/p}\\&\quad\le\sum_{m=0}^{\infty}\frac{\delta^{qm}}{m!}\left(\frac{b_n}{n}\right)^{m}a\left(\frac{n}{b_n}\right)^{md/p}(\mathbb{E}L^{m}(n,0))^{1/p},\end{aligned}$$

by (1.31),

$$\begin{aligned}&\limsup_{n\to\infty}\frac{1}{b_n}\log\sum_{m=0}^{\infty}\frac{1}{(2m)!}\left(\frac{b_n}{n}\right)^{2m}a\left(\frac{n}{b_n}\right)^{2md/p}(\mathbb{E}L^{2m}(n,0))^{1/p}\\&\quad\le(2\pi)^{-\alpha d/(\alpha p-d)}\rho_1^{\alpha/(\alpha p-d)}\delta^{\alpha pq/(\alpha p-d)}.\end{aligned}$$

Notice that we can make the right-hand side arbitrarily small by controlling $\delta$. To prove (6.7), therefore, we need only to show that, for any $\theta > 0$,

$$\lim_{n\to\infty}\frac{1}{b_n}\log\sum_{m=0}^{\infty}\frac{\theta^{2m}}{(2m)!}\left(\frac{b_n}{n}\right)^{2m}a\left(\frac{n}{b_n}\right)^{2md/p}(\mathbb{E}\tilde{\xi}^{2m}(n,0))^{1/p}=0. \tag{6.8}$$

For any $A\subset(\mathbb{Z}^+)^p$, write

$$\tilde{\xi}(A)=\sum_{(k_1,\dots,k_p)\in A}(\varepsilon^1_{k_1}\tau^1_k)\cdots(\varepsilon^p_{k_p}\tau^p_k)\mathbb{1}_{\{\tilde{S}_1(k_1)+\cdots+\tilde{S}_p(k_p)=0\}}.$$

We claim that, for any $m \ge 1$ and any $A \subset B$,

$$\mathbb{E}\tilde{\xi}^{2m}(A) \le \mathbb{E}\tilde{\xi}^{2m}(B). \tag{6.9}$$

Indeed, using the Jensen inequality conditionally on

$$\{(\varepsilon_{k_1}\tau_{k_1}, \dots, \varepsilon_{k_p}\tau_{k_p}); (k_1, \dots, k_p) \in B \setminus A\}$$

gives

$$\begin{aligned}
&\mathbb{E}\tilde{\xi}^{2m}(B) \\
&\quad \ge \mathbb{E}\Bigg[\sum_{k_1,\dots,k_p \in A} (\varepsilon^1_{k_1}\tau^1_k) \cdots (\varepsilon^p_{k_p}\tau^p_k) \mathbb{1}_{\{\tilde{S}_1(k_1)+\cdots+\tilde{S}_p(k_p)=0\}} \\
&\qquad\quad + \sum_{(k_1,\dots,k_p)\in B\setminus A} \mathbb{E}[(\varepsilon^1_{k_1}\tau^1_k) \cdots (\varepsilon^p_{k_p}\tau^p_k)] \mathbb{1}_{\{\tilde{S}_1(k_1)+\cdots+\tilde{S}_p(k_p)=0\}}\Bigg]^{2m} \\
&\quad = \mathbb{E}\tilde{\xi}^{2m}(A).
\end{aligned}$$

In particular, by (6.9), we have $\mathbb{E}\tilde{\xi}^{2m}(n,0) \le \mathbb{E}\tilde{\xi}^{2m}(n+1,0)$ for all $n \ge 1$. Write $r_n = [nb_n^{-1}] + 1$. By Lemma 4.5,

$$\begin{aligned}
&\sum_{m=0}^{\infty} \frac{\theta^{2m}}{(2m)!}\left(\frac{b_n}{n}\right)^{2m} a\left(\frac{n}{b_n}\right)^{2md/p} (\mathbb{E}\tilde{\xi}^{2m}(n,0))^{1/p} \\
&\qquad \le \left\{\sum_{m=0}^{\infty} \frac{\theta^{2m}}{(2m)!}\left(\frac{b_n}{n}\right)^{2m} a\left(\frac{n}{b_n}\right)^{2md/p} (\mathbb{E}\tilde{\xi}^{2m}(r_n,0))^{1/p}\right\}^{[b_n]+1}.
\end{aligned} \tag{6.10}$$

Notice that

$$\begin{aligned}
\mathbb{E}\tilde{\xi}^2(r_n,0) &= \mathbb{E}\sum_{k_1,\dots,k_p=0}^{r_n} \mathbb{1}_{\{\tilde{S}_1(k_1)+\cdots+\tilde{S}_p(k_p)=0\}} \\
&= O\left\{\left(\frac{n}{b_n}\right)^p a\left(\frac{n}{b_n}\right)^{-d}\right\} \qquad (n \to \infty).
\end{aligned} \tag{6.11}$$

In addition,

$$\begin{aligned}
|\tilde{\xi}(r_n,0)| &\le \sum_{k_1,\dots,k_p=0}^{r_n} \tau^1_{k_1} \cdots \tau^p_{k_p} \mathbb{1}_{\{S_1(2k_1)+\cdots+S_p(2k_p)=0\}} \\
&\overset{d}{=} \sum_{k_1,\dots,k_p=0}^{r_n} \tau^1_{2k_1} \cdots \tau^p_{2k_p} \mathbb{1}_{\{S_1(2k_1)+\cdots+S_p(2k_p)=0\}} \\
&\le \sum_{k_1,\dots,k_p=0}^{r_n} \tau^1_{k_1} \cdots \tau^p_{k_p} \mathbb{1}_{\{S_1(k_1)+\cdots+S_p(k_p)=0\}} \le L(T^*_{r_n+1}, 0),
\end{aligned}$$

where $T_n^* = \max_{1\le j\le p} T_n^j$ and the where last step follows from (4.39).

Using (4.31) in Lemma 4.5 conditionally on $T_n^*$, we have that for any $\delta > 0$ there is a $C > 0$ such that, as $n$ is sufficiently large,

$$\begin{aligned}
&\mathbb{E}L^{2m}(T^*_{r_n+1}, 0)\\
&\quad \le C[(2m)!]^p \delta^{2m}\mathbb{E}\{(T_n^*)^p a(T_n^*)^{-d}\}^{2m}\\
&\quad \le C[(2m)!]^p (\lambda\delta)^{2m}\left(\frac{n}{b_n}\right)^{2mp} a\left(\frac{n}{b_n}\right)^{-2md}\\
&\qquad + C[(2m)!]^p \delta^{2m}\mathbb{E}\{(T_n^*)^p a(T_n^*)^{-d}\}^{2m}\mathbb{1}_{\{T_{r_n+1}\ge 2r_n\}}
\end{aligned}$$

for any $m \ge 1$, where $\lambda$ can be any fixed number greater than $2^{(\alpha p-d)/\alpha}$, and the second step follows partially from the regularity given in (1.4).

Using (6.1) to control the second term on the right-hand side, we obtain

$$\mathbb{E}\tilde{\xi}^{2m}(r_n, 0) \le C'[(2m)!]^p (\lambda\delta)^{2m}\left(\frac{n}{b_n}\right)^{2mp} a\left(\frac{n}{b_n}\right)^{-2md}$$

for all $m \ge 1$ as $n$ is sufficiently large.

We now take $\delta$ small enough so $\lambda\delta < \theta^{-1}$. By the dominated convergence theorem and by (6.11),

$$\sum_{m=0}^{\infty} \frac{\theta^{2m}}{(2m)!}\left(\frac{b_n}{n}\right)^{2m} a\left(\frac{n}{b_n}\right)^{2md/p} (\mathbb{E}\tilde{\xi}^{2m}(r_n, 0))^{1/p} \longrightarrow 1 \qquad (n \to \infty).$$

Thus, (6.8) follows from (6.10).

## APPENDIX

Let $f \in \mathcal{L}^1(\mathbb{R}^d)$ be a symmetric and nonnegative function on $\mathbb{R}^d$ and write

$$\bar{f}(x) = \int_{\mathbb{R}^d} e^{i\lambda\cdot x} f(\lambda)\, d\lambda, \qquad x \in \mathbb{R}^d.$$

Let $\mathcal{F}_\Psi$ be defined as in (2.4) and write

$$\rho(f) = \sup_{\|g\|_2=1} \iint_{\mathbb{R}^d\times\mathbb{R}^d} f(\lambda-\gamma)\frac{g(\lambda)g(\gamma)}{\sqrt{1+\Psi(\lambda)}\sqrt{1+\Psi(\gamma)}}\, d\lambda\, d\gamma.$$

In this section we prove the following lemma which was used in Section 2.

LEMMA A.1.

$$\sup_{g\in\mathcal{F}_\Psi}\left\{\frac{1}{\rho(f)}\int_{\mathbb{R}^d}\bar{f}(x)g^2(x)\,dx - \int_{\mathbb{R}^d}|\widehat{g}(\lambda)|^2\Psi(\lambda)\,d\lambda\right\} = 1, \tag{A.1}$$

*where*

$$\widehat{g}(\lambda) = \int_{\mathbb{R}^d} g(x)e^{i\lambda\cdot x}\,dx, \qquad \lambda \in \mathbb{R}^d.$$

PROOF. It is easy to see that the supremum in the definition of $\rho(f)$ can be taken only over the symmetric functions $g$. Write $Q(\lambda) = (1 + \Psi(\lambda))^{-1}$ and $h(\lambda) = g(\lambda)/\sqrt{Q(\lambda)}$. Then

$$\|h\|^2_{\mathcal{L}^2(Q)} \equiv \int_{\mathbb{R}^d} f^2(\lambda) Q(\lambda)\, d\lambda = 1,$$

$$\rho(f) = \sup_{\|h\|_{\mathcal{L}^2(Q)}=1} \int\int_{\mathbb{R}^d\times\mathbb{R}^d} f(\lambda-\gamma) Q(\lambda) h(\lambda) Q(\gamma) h(\gamma)\, d\lambda\, d\gamma \tag{A.2}$$
$$= \sup_{\|h\|_{\mathcal{L}^2(Q)}=1} \int_{\mathbb{R}^d} f(\lambda) U(\lambda)\, d\lambda,$$

where

$$U(\lambda) = \int_{\mathbb{R}^d} Q(\lambda+\gamma) h(\lambda+\gamma) Q(\gamma) h(\gamma)\, d\gamma.$$

Notice that

$$\begin{aligned} V(x) &\equiv \frac{1}{(2\pi)^d} \int_{\mathbb{R}^d} U(\lambda) e^{-i\lambda\cdot x}\, d\lambda \\ &= \frac{1}{(2\pi)^d} \int_{\mathbb{R}^d} e^{-i\lambda\cdot x} d\lambda \int_{\mathbb{R}^d} Q(\gamma+\lambda) Q(\gamma) h(\gamma+\lambda) h(\gamma)\, d\gamma \\ &= \frac{1}{(2\pi)^d} \int\int_{\mathbb{R}^d\times\mathbb{R}^d} e^{-i(\lambda-\gamma)\cdot x} Q(\lambda) h(\lambda) Q(\gamma) h(\gamma)\, d\lambda\, d\gamma \\ &= \frac{1}{(2\pi)^d} \left| \int_{\mathbb{R}^d} e^{ix\cdot\gamma} Q(\gamma) h(\gamma) d\gamma \right|^2 \end{aligned} \tag{A.3}$$

and that $Q(\lambda)$ is the Fourier transform of the killed Green's function

$$G(x) = \int_0^\infty e^{-t} p_t(x)\, dt, \qquad x \in \mathbb{R}^d,$$

where $p_t$ is the density function of $Y_t$. Under the substitution

$$g(x) = \frac{1}{(2\pi)^d} \int_{\mathbb{R}^d} e^{-i\lambda\cdot x} h(\lambda)\, d\lambda,$$

we have

$$\int_{\mathbb{R}^d} e^{ix\cdot\gamma} Q(\gamma) h(\gamma)\, d\gamma = (2\pi)^d \int_{\mathbb{R}^d} G(y-x) h(y)\, dy = (2\pi)^d G h(x) \tag{A.4}$$

and, therefore,

$$\|h\|^2_{\mathcal{L}^2(Q)} = (2\pi)^d \int_{\mathbb{R}^d\times\mathbb{R}^d} G(x-y) g(x) g(y)\, dx\, dy = (2\pi)^d \langle g, Gg\rangle. \tag{A.5}$$

By (A.2)–(A.5), and by Parseval's identity,

$$\begin{aligned}\rho(f) &= \sup_{\langle g,Gg\rangle=(2\pi)^{-d}} \int_{\mathbb{R}^d} \bar{f}(x)V(x)\,dx \\ &= (2\pi)^d \sup_{\langle g,Gg\rangle=(2\pi)^{-d}} \int_{\mathbb{R}^d} \bar{f}(x)[Gg(x)]^2\,dx.\end{aligned} \tag{A.6}$$

Write $h(x)=Gg(x)$ and recall the resolvent identity

$$I = G - \mathcal{A}\circ G,$$

where $I$ is the identical operator and where $\mathcal{A}$ is the infinitesimal generator of the Markov process $Y(t)$. It is a well-known fact (page 24, [4]) that the linear operator $\mathcal{A}$ is determined by

$$\mathcal{A}h(x) = -\int_{\mathbb{R}^d} \Psi(\lambda)\widehat{h}(\lambda)e^{-i\lambda\cdot x}\,d\lambda, \qquad x\in\mathbb{R}^d.$$

Hence,

$$\begin{aligned}\langle g,Gg\rangle = \langle h-\mathcal{A}h,h\rangle &= \|h\|_2 + \int_{\mathbb{R}^d}\psi(\lambda)|\widehat{h}(\lambda)|^2\,d\lambda \\ &= \|h\|_2 + \|\widehat{h}\|^2_{\mathcal{L}^2(\Psi)},\end{aligned}$$

where

$$\|h\|^2_{\mathcal{L}^2(\Psi)} = \int_{\mathbb{R}^d}\Psi(\lambda)|h(\lambda)|^2\,d\lambda.$$

From (A.6), we have

$$\rho(f) = (2\pi)^d \sup_{\|h\|^2+\|\widehat{h}\|^2_{\mathcal{L}^2(\Psi)}=(2\pi)^{-d}} \int_{\mathbb{R}^d}\bar{f}(x)h^2(x)\,dx. \tag{A.7}$$

In addition, it is easy to see that the function

$$M_f(\theta) = \sup_{g\in\mathcal{F}_\Psi}\left\{\theta\int_{\mathbb{R}^d}\bar{f}(x)g^2(x)\,dx - \int_{\mathbb{R}^d}|\widehat{g}(\lambda)|^2\Psi(\lambda)\,d\lambda\right\}, \qquad \theta>0,$$

is nonnegative, nondecreasing and continuous on $(0,\infty)$. By (A.7), given $0<\varepsilon<\rho(f)$, there is $h_o$ such that $\|h_o\|^2+\|\widehat{h_o}\|^2_{\mathcal{L}^2(\psi)}=(2\pi)^{-d}$ and that

$$\int_{\mathbb{R}^d}\bar{f}(x)h_o^2(x)\,dx > (2\pi)^{-d}(\rho(f)-\varepsilon).$$

Hence,

$$\begin{aligned}M_f(\theta)\left(\frac{1}{\rho(f)-\varepsilon}\right) &\ge \frac{1/(\rho(f)-\varepsilon)\int_{\mathbb{R}^d}\bar{f}(x)h_o^2(x)\,dx - \int_{\mathbb{R}^d}\Psi(\lambda)|\widehat{h_o}(\lambda)|^2\,d\lambda}{\int_{\mathbb{R}^d}|h_o(x)|^2\,dx} \\ &= \frac{(2\pi)^{-d} - \int_{\mathbb{R}^d}\Psi(\lambda)|\widehat{h_o}(\lambda)|^2\,d\lambda}{\int_{\mathbb{R}^d}|h_o(x)|^2\,dx} = 1.\end{aligned}$$

Letting $\varepsilon \to 0^+$ on the left-hand side gives that

$$M_f\left(\frac{1}{\rho(f)}\right) \ge 1.$$

On the other hand, by (A.7),

$$\begin{aligned} M_f\left(\frac{1}{\rho(f)}\right) &= \sup_{g\in\mathcal{F}_d}\left\{\frac{1}{\rho(f)}\int_{\mathbb{R}^d}\bar f(x)g^2(x)\,dx - \int_{\mathbb{R}^d}\Psi(\lambda)|\widehat{g}(\lambda)|^2\,d\lambda\right\} \\ &\le \sup_{g\in\mathcal{F}_d}\left\{\frac{1}{\rho(f)}\rho(f)\left[1+\int_{\mathbb{R}^d}\Psi(\lambda)|\widehat{g}(\lambda)|^2\,d\lambda\right] - \int_{\mathbb{R}^d}\Psi(\lambda)|\widehat{g}(\lambda)|^2\,d\lambda\right\} = 1. \end{aligned}$$

□

**Acknowledgment.** The author is grateful to an anonymous referee for his/her careful reading, as well as for pointing out some misprints and English errors in the first draft of this paper.

Department of Mathematics
University of Tennessee
Knoxville, Tennessee 37996-1300
USA
E-mail: xchen@math.utk.edu
URL: http://www.math.utk.edu/~xchen